\documentclass[a4paper]{article}
\usepackage{amsthm,amsmath,amssymb,mathrsfs}

\usepackage{color}
\usepackage{cite}

\allowdisplaybreaks[4]
\newtheorem{teor}{Theorem}[section]
\newtheorem{defin}[teor]{Definition}
\newtheorem{lemm}[teor]{Lemma}
\newtheorem{osse}[teor]{Remark}
\newtheorem{prop}[teor]{Proposition}
\newtheorem{defi}[teor]{Definition}
\newtheorem{coro}[teor]{Corollary}
\newtheorem{prob}[teor]{Problem}

\newcommand{\bele}{\begin{lemm}\begin{sl}}
\newcommand{\enle}{\end{sl}\end{lemm}}
\newcommand{\bedef}{\begin{defi}\begin{sl}}
\newcommand{\eddef}{\end{sl}\end{defi}}
\newcommand{\bete}{\begin{teor}\begin{sl}}
\newcommand{\ente}{\end{sl}\end{teor}}
\newcommand{\beos}{\begin{osse}\begin{rm}}
\newcommand{\eddos}{\end{rm}\end{osse}}
\newcommand{\bepr}{\begin{prop}\begin{sl}}
\newcommand{\empr}{\end{sl}\end{prop}}
\newcommand{\bepro}{\begin{prob}\begin{rm}}
\newcommand{\empro}{\end{rm}\end{prob}}
\newcommand{\bede}{\begin{defin}\begin{sl}}
\newcommand{\edde}{\end{sl}\end{defin}}
\newcommand{\beco}{\begin{coro}\begin{sl}}
\newcommand{\enco}{\end{sl}\end{coro}}

\newcommand{\quext}{\quad\text}

\newcommand{\de}{\partial}

\newcommand{\RR}{\mathbb{R}}

\newcommand{\NN}{\mathbb{N}}

\newcommand{\beeq}[1]{\begin{equation}\label{#1}}
\newcommand{\eddeq}{\end{equation}}

\newcommand{\beeqa}[1]{\begin{eqnarray}\label{#1}}
\newcommand{\eddeqa}{\end{eqnarray}}

\newcommand{\beal}[1]{\begin{align}\label{#1}}
\newcommand{\eddal}{\end{align}}

\newcommand{\bespl}[1]{\begin{split}\label{#1}}
\newcommand{\edspl}{\end{split}}

\newcommand{\bega}[1]{\begin{gather}\label{#1}}
\newcommand{\edga}{\end{gather}}

\newcommand{\beeqax}{\begin{eqnarray*}}
\newcommand{\eddeqax}{\end{eqnarray*}}

\newcommand{\beeqao}{\begin{eqnarray}\no}
\newcommand{\bealo}{\begin{align}\no}
\newcommand{\besplo}{\begin{split}\no}
\newcommand{\begao}{\begin{gather}\no}

\def\qed{\ifmmode % if math mode, assume display: omit penalty etc.
  \else \leavevmode\unskip\penalty9999 \hbox{}\nobreak\hfill
  \fi
  \quad\hbox{\hskip.5em\vrule width.4em height.6em depth.05em\hskip.1em}}
\def\endproofsym{\qed}
\renewenvironment{proof}[1][Proof]{\trivlist\item[\hskip\labelsep{\hskip0pt
    {\normalfont\scshape#1.}\hskip .321429\parindent}]\ignorespaces}
{\endproofsym\endtrivlist}
\def\endnobox{\def\endproofsym{}\end{proof}\def\endproofsym{\qed}}

\newcommand{\no}{\nonumber}

\newcommand{\duav}[1]{\langle{#1}\rangle}

\newcommand{\perogni}{\forall\,}
\newcommand{\esiste}{\exists\,}
\newcommand{\itt}{\int_0^t}

\newcommand{\io}{\int_\Omega}
\newcommand{\ibaro}{\int_{\overline{\Omega}}}

\newcommand{\iTT}{\int_0^T}

\newcommand{\iTo}{\iTT\!\io}

\newcommand{\epsi}{\varepsilon}

\def\R{\mathbb R}

\newcommand{\bn}{\boldsymbol{n}}

\newcommand{\dn}{\partial_{\bn}}
\newcommand{\fhi}{\varphi}

\newcommand{\lhs}{left hand side}
\newcommand{\rhs}{right hand side}

\DeclareMathOperator{\dive}{div}
\DeclareMathOperator{\deriv}{d}

\DeclareMathOperator{\dist}{dist}

\newcommand{\LDV}{L^2(0,T;V)}

\let\TeXchi\chi
\def\chi{{\setbox0 \hbox{\mathsurround0pt
$\TeXchi$}\hbox{\raise\dp0 \copy0 }}}

\newcommand{\zzn}{_{0,n}}

\newcommand{\calX}{{\mathcal X}}

\newcommand{\calA}{{\mathcal A}}

\newcommand{\calE}{{\mathcal E}}

\newcommand{\calM}{{\mathcal M}}

\newcommand{\calR}{{\mathcal R}}
\newcommand{\calB}{{\mathcal B}}

\newcommand{\zzu}{_{0,1}}
\newcommand{\zzd}{_{0,2}}

\newcommand{\barO}{\overline{\Omega}}

\newcommand{\dit}{\deriv\!t}
\newcommand{\dis}{\deriv\!s}
\newcommand{\dix}{\deriv\!x}

\newcommand{\dir}{\deriv\!r}

\newcommand{\ddt}{\frac{\deriv\!{}}{\dit}}
\newcommand{\ddr}{\frac{\deriv\!{}}{\dir}}

\newenvironment{bettirev}{\color{blue}}{\color{black}}
\newcommand{\bber}{\begin{bettirev}}
\newcommand{\eber}{\end{bettirev}}

\newenvironment{michelarev}{\color{red}}{\color{black}}
\newcommand{\III}{\begin{michelarev}}
\newcommand{\EEE}{\end{michelarev}}

\numberwithin{equation}{section}
\begin{document}

\title{On a Class of Sixth-order Cahn--Hilliard Type Equations \\
with Logarithmic Potential}

\author{%
Giulio Schimperna\\
Dipartimento di Matematica, Universit\`a di Pavia,\\
Via Ferrata~1, 27100 Pavia, Italy\\
E-mail: {\tt giusch04@unipv.it}\\
\and
Hao Wu \\
School of Mathematical Sciences\\
Shanghai Key Laboratory for Contemporary Applied Mathematics\\
Fudan University\\
Han Dan Road 220, Shanghai 200433, China\\
E-mail: {\tt haowufd@fudan.edu.cn}
}

%\date{}

\maketitle
\begin{abstract}
We consider a class of six-order Cahn--Hilliard equations with logarithmic type potential. This system is closely connected with some important phase-field models 
relevant in different applications, for instance, the functionalized Cahn--Hilliard equation that describes phase separation
in mixtures of amphiphilic molecules in solvent, and the Willmore regularization of Cahn--Hilliard equation for anisotropic crystal and epitaxial growth. 
The singularity of the configuration potential guarantees that the solution always stays in the physical relevant domain $[-1,1]$. Meanwhile, the resulting 
system is characterized by some highly singular diffusion terms that make the mathematical analysis more involved. We prove existence and uniqueness of 
global weak solutions and show their parabolic regularization property for any  positive time. Besides, we investigate long-time behavior of the system, 
proving existence of the global attractor for the associated
dynamical process in a suitable complete metric space.
\end{abstract}


\vspace{2mm}

\noindent {\bf AMS (MOS) subject clas\-si\-fi\-ca\-tion:}~~35K35, 35K55, 35A01, 47H05.

%%%%%%%%%%%%%%%%%%%%%%%%%%%%%%%%%%%%%%%%%%%%%%%%%%%%%%%%%%%%%%%%%%%%%%%%%%%%%%%%%%%%%%%%%%%%%
%
% INTRODUCTION
%
%%%%%%%%%%%%%%%%%%%%%%%%%%%%%%%%%%%%%%%%%%%%%%%%%%%%%%%%%%%%%%%%%%%%%%%%%%%%%%%%%%%%%%%%%%%%%

\section{Introduction}
\label{sec:intro}
\setcounter{equation}{0}

In this paper, we study the following class of parabolic systems:
\begin{align}
 \label{eq:mu}
  & \partial_t u  =  \Delta \mu &\quext{in }\,(0,T)\times\Omega,\\
 \label{eq:u}
  & \mu =  - \Delta \omega + f'(u) \omega  + \eta \omega  &\quext{in }\,(0,T)\times\Omega,\\
 \label{eq:om}
  & \omega = - \Delta u + f(u) &\quext{in }\,(0,T)\times\Omega,
\end{align}
subject to the boundary and initial conditions
\begin{align}
   & \dn u=\dn \Delta u=\dn \mu=0 &\text{on }\,(0,T)\times \partial \Omega,\label{bd}\\
   & u|_{t=0}=u_0(x) &\text{in }\,\Omega.\label{ini}
\end{align}
Here, $\Omega$ is assumed to be a smooth bounded domain in~$\RR^3$ with boundary~$\partial \Omega$ and $T>0$ is a given final time of arbitrary magnitude. The restriction 
to the three-dimensional setting is motivated by physical applications. Actually, similar, or even better results are expected to hold in space dimensions one and two. 
We denote the outward unit normal vector to $\partial \Omega$ by $\bn$, and by $\dn$ the outward normal derivative. 
The homogeneous Neumann boundary conditions in \eqref{bd}
are standardly used in the frame of Cahn--Hilliard models due to their mathematical simplicity. Besides, they are related to some basic features of the problem, e.g., 
the variational structure and the mass conservation. Other types of boundary conditions for $u$ might be considered as well, 
for instance, the periodic boundary conditions for a cubical domain like $\Omega=[0,1]^3$.
We also observe that \eqref{bd} is equivalent to asking 
$\dn u=\dn \omega=\dn \mu=0$ on $(0,T)\times \partial \Omega$. 

System \eqref{eq:mu}--\eqref{eq:om} is a variant of the classical fourth-order Cahn--Hilliard equation \cite{CH} that describes possible separation process of binary mixtures. The variable $u$ has the meaning of an order parameter representing the local proportion of one of the two components of the binary material.
For simplicity, we may assume $u$ to be normalized in such a way that the pure states
correspond to the values  $u = \pm 1$, whereas $-1<u<1$ denotes the (local)
presence of a mixture. The function $\mu$ in \eqref{eq:u} is an auxiliary variable generally termed as ``chemical potential'', which is the first variational derivative of the following energy functional
\begin{equation}\label{defiE}
  \calE(u) = \io \bigg[ \frac12 \big| - \Delta u + f(u) \big|^2 + \eta \Big( \frac12 |\nabla u|^2 + F(u) \Big) \bigg] \dix.
\end{equation}
The energy $\calE(u)$ can be seen as an extension of the Cahn--Hilliard free energy (cf.~\cite{CH})
\begin{align}
\calE_{CH}(u)= \io \frac12 |\nabla u|^2 + F(u) \dix,\no
\end{align} 
whose first variational derivative is denoted by $\omega$ as in \eqref{eq:om}. 
$F$ is a configuration potential function that usually has a double-well structure, with its derivative denoted by $f=F'$.

Throughout the paper, the parameter $\eta$ is assumed to be a constant for simplicity, and it can take values in $\RR$. We note that the value/sign of $\eta$ is important in modeling and application. When $\eta=0$, $\calE$ reduces to the well known Willmore functional in the phase-field formulation that approximates the Canham--Helfrich bending energy of surfaces. This was efficiently used to study deformations of elastic vesicles subject to possible volume/surface constraints \cite{DLW04,DLRW05}. Next, for $\eta>0$, $\calE$ is related to the Willmore regularization of the Cahn--Hilliard energy $\calE_{CH}$, which was introduced for instance, in \cite{CS,TLVW09} to investigate strong anisotropy effects (corresponding to some non-constant coefficient $\eta$) arising during the growth and coarsening of thin films. When $\eta<0$, $\calE$ is referred to as the functionalized Cahn--Hilliard (FCH) free energy, which was derived from models for mixtures with an amphiphilic structure \cite{GS90} and for nanoscale morphology changes in functionalized polymer chains \cite{PW09}. In this case, $\calE$ reflects the balance between the square of the variational derivative of $\calE_{CH}$ against itself, such that it can incorporate the propensity of the amphiphilic surfactant phase to drive the creation of interfaces and naturally produce stable bilayers, or homoclinic interfaces with an intrinsic width (cf. \cite{DP13}). Minimization problems, bilayer structures, pearled patterns, and network bifurcations related to the FCH energy have been extensively studied in \cite{DP15,PZ13,PW17} and the references therein.

In this paper, we shall analyze the evolution problem \eqref{eq:mu}--\eqref{ini} associated with the energy functional $\calE$ given by \eqref{defiE}, for arbitrary $\eta\in \RR$ in a uniform way. The problem admits a natural variational
interpretation such that it can be restated as
\begin{align}\label{gf:1}
  & \partial_t u = \dive (M(u) \nabla \mu),\quad \text{with}\quad 
   \mu = \frac{\delta \calE}{\delta u},
\end{align}
where $M(u)$ is the mobility function (in our present case, $M\equiv 1$). 
The first relation in \eqref{gf:1} represents a continuity equation that corresponds to the conservation of mass. Namely,
the mass flux $M(u)\nabla \mu$ is proportional to the gradient of the
chemical potential $\mu$ through the mobility $M$. As a consequence of the no-flux boundary condition~\eqref{bd},
integrating over $\Omega$, we obtain 
\begin{equation}\label{st:mean}
  \ddt \io u \dix = 0.
\end{equation}
This property is consistent with the physical ansatz that the total amount
of each component of the material is conserved in time. Besides, formally
testing the first of \eqref{gf:1} by $\mu$ and using again the boundary condition \eqref{bd}, we can derive the following energy dissipation law 
\begin{equation}\label{diss:E}
  \ddt \calE(u) + \io M(u) | \nabla \mu |^2 \dix = 0.
\end{equation}
In particular, for our case $M\equiv 1$, \eqref{eq:mu} turns out to be the usual $(H^{1})'$-gradient flow of $\calE$. 

In order to understand the main features of our system and the
mathematical difficulties occurring in its analysis,
we now specify our assumptions on the potential function $F$ and its derivative $f$. In view of the physical interpretation of $u$, only the values $u \in [-1, 1]$ are admissible and this constraint should be somehow enforced in the equations. Due to the lack of maximum principle for higher order equations, a common way to do this
consists in taking $f$ as the derivative of a non-convex configuration
potential $F$ of logarithmic type (sometimes referred to as the Flory--Huggins potential). Without loss of generality, in this paper we set
\begin{equation}\label{defiF}
  F(r) = \frac12 (1 + r) \ln  (1 + r) + \frac12 (1 - r) \ln (1 - r) - \frac\lambda2 r^2,
   \quad r\in(-1,1),
\end{equation}
where $\lambda \in \mathbb{R}$ is a constant. Then its derivative $f$ is given by
\begin{equation}\label{defif}
  f(r) = F'(r) = \frac12 \ln (1 + r) - \frac12 \ln (1 - r) - \lambda r,\quad r\in(-1,1).
\end{equation}
It is worth noting that, when the (given) constant $\lambda$ is large positive, $F$ may be non-convex.
This is physically reasonable in view of the fact that the minima of $F$
correspond to the energetically more favorable configurations attained
in proximity of the pure states $u\sim \pm 1$. In the terminology of convex analysis, $F$ is a {\sl $\lambda$-convex}\/
function (i.e., convex if the perturbation induced by the
$\lambda$-term is neglected) whose effective domain is the interval $[-1,1]$.
Indeed, $F(r)$ can be extended by continuity to $r=\pm 1$ and thought to
be identically $+\infty$ for $r\not \in [-1,1]$. In practice, the singular potential $F$ is often approximated by regular ones of polynomial, for instance, $F(r)=\frac{1}{4}(r^2-1)^2$. 

There is an increasing interest on the study of six-order Cahn--Hilliard type equations appearing as models of various physical phenomena. For instance, focusing on their mathematical analysis, we may refer to \cite{LZ18,PZa11,PZa13,SP0} for dynamics of oil-water-surfactant mixtures, to \cite{KNR12,KR12} for the faceting of growing surfaces, and to \cite{GW14,GW15,Mir15,Mir16,WW10} for the phase-field-crystal equation modeling crystal growth. 
Concerning our problem \eqref{eq:mu}--\eqref{ini}, although there have been extensive numerical studies for the model \cite{CLSWW,CJPWW,FGLWW,GXX}, to the best of our knowledge, only few works 
from the theoretical point of view have been done in the literature.
For $\eta=1$, in \cite{Mir1}, the author considers the case of a regular (polynomial) potential $F$ with constant mobility, proving existence and uniqueness of global weak solutions, existence of the finite dimensional global attractor and exponential attractors, while existence of weak solutions with a degenerate phase-dependent mobility in the two dimensional case is obtained in \cite{DCZ19}. For the functionalized Cahn--Hilliard equation (i.e., $\eta<0$) subject to periodic boundary conditions, in \cite{DLP19}, the authors proved existence of global weak solutions in the case of regular potential and degenerate mobility. Besides, for $\eta\in \R$, $M=1$ and a regular potential, existence and uniqueness of global solutions in the Gevrey class were established in \cite{CWWY}, again in the periodic setting. We note that in those contributions \cite{CWWY,DCZ19,DLP19,Mir1} mentioned above, the potential $F$ is always assumed to be a regular one. Indeed, this choice plays a crucial role in the mathematical analysis therein.   

Our aim in this paper is twofold. More precisely, for problem \eqref{eq:mu}--\eqref{ini} with the logarithmic potential $F$ given by \eqref{defiF}, we shall prove for $\eta, \lambda\in \RR$:
\begin{itemize}
\item[(A)] existence and uniqueness of global weak solutions in a suitable regularity class (Theorem \ref{thm:main}) and the parabolic regularization property of solutions for any $t>0$ (Theorem \ref{thm:furth});
\item[(B)] long-time behavior of the dynamical system associated with problem \eqref{eq:mu}--\eqref{ini}, in terms of existence of the global attractor in a suitable complete metric space (Theorem \ref{thm:attra}). 
\end{itemize}
The main difficulty in mathematical analysis is due to the singular character of $F$ and its interaction with higher-order derivatives in the present six-order equation. We note that although a rather complete characterization on well-posedness, regularity as well as long-time behavior of the classical fourth-order Cahn--Hilliard equation with singular (logarithmic) potential can be found in the literature (see e.g., \cite{AM,CMZ,MZ} and references therein), the situation is far less understood in the six-order system. In this direction, we refer to \cite{Mir15, Mir16} for possibly the simplest case related to the phase-field-crystal model and to \cite{SP0} for the case with further nonlinear (but regular) diffusion, modeling oil-water-surfactant mixtures. In our case, the situation is much more involved due to the presence of the singular diffusion term $-\Delta f(u)$ in the chemical potential $\mu$ 
(coming from the combination of \eqref{eq:u} with \eqref{eq:om}, see
\eqref{eq:uom} below). This is similar to the situation faced in the
paper \cite{SP} for the fourth-order Cahn--Hilliard equation, where an interpretation of the singular diffusion terms by means of variational inequalities and duality methods was used in order to overcome the lack of
available regularity properties. 

Hence, several techniques have to be exploited to handle those mathematical difficulties due to the highly singular and nonlinear structure of the system. 

First, in order to prove the existence of weak solutions, we take advantage of the six-order Laplacian and introduce a nonstandard ``singular" approximation for the original problem \eqref{eq:mu}--\eqref{ini}, using a sequence of singular functionals induced by indicator functions with domain $[-1+1/n,\,1-1/n]$ (see \eqref{eq:mu:d}--\eqref{eq:eta:d}). This enables us to avoid the singularities due to $f(u)$ by paying the price to deal with an additional singular term $\zeta_n$, which is indeed simpler since it acts only 
on the solution $u$ and does not involve any derivative of it. Then we can apply a similar argument in \cite{SP0} to conclude the existence of approximate solutions. On the other hand, in view of the specific variational structure
of the system, we are able to deduce sufficient \textit{a priori} estimates (see Section \ref{sec:apriori}) and
show that the singular terms satisfy a sort of equiintegrability
property in such a way that they can be interpreted in the classical sense
(i.e., pointwise), avoiding the use of variational inequalities and
duality methods. The approximating scheme turns out to be consistent with those formal \textit{a priori} estimates and this allows us to pass to the limit $n \nearrow +\infty$ to obtain the existence of a weak solution. 

Second, in order to prove the uniqueness, we define a new class of solutions based on abstract subdifferential methods (and, hence, called ``subdifferential
solutions'') and prove that this class is wider than that of the
``weak solutions'' considered before. Then, uniqueness of subdifferential
solutions can be shown by means of chain rule formulas and convex
analysis methods, and it immediately implies uniqueness of weak solutions.

Third, after showing that the global weak solutions satisfy parabolic time-regularization properties for any $t>0$, we can proceed to study their long-time behavior. Namely, we prove existence of the global attractor for the dynamical system generated by solution trajectories in the natural phase space, which is chosen precisely as the space for the initial data satisfying the conditions required for having existence of a weak solution. Here, we recall that the global attractor is the smallest compact set in a certain phase space that is invariant under the semiflow $S(t)$ generated by the solution $u(t)$ and attracts all bounded sets of initial data as time goes to infinity \cite{BV,Tem}. It is worth mentioning that, as in similar
situations involving singular nonlinearities \cite{RS} (see also \cite{WZ07}), the phase space we consider here has no Banach structure, but is just a complete metric spaces with a proper distance also acting on the singular term $f(u)$. In this framework, we are able to prove the existence of a compact absorbing set with respect to the metric of the phase space, which finally yields the existence of the compact attractor.

The plan of the paper is as follows: Section \ref{sec:main} is devoted to presenting
the weak formulation of our problem \eqref{eq:mu}--\eqref{ini} and stating our main results regarding well-posedness, parabolic regularization of solutions, as well as long-time behavior.
Their proofs are postponed to the sequel. In particular, in Section~\ref{sec:apriori} we derive a sequence of \textit{a priori} estimates
required for proving the existence of solutions. This is, indeed, the key point of our argument where several mathematical difficulties are concentrated. Then, the proof of existence will be completed in Section~\ref{sec:exi} by using a suitable approximation scheme and asymptotic compactness methods. Uniqueness will be discussed in Section~\ref{sec:uniq} by restating our problem as an abstract evolution equation ruled by a subdifferential operator and applying some convex analysis tools. Finally, regularization properties of solutions and their long-time behavior will be analyzed in Section~\ref{sec:long}.

%%%%%%%%%%%%%%%%%%%%%%%%%%%%%%%%%%%%%%%%%%%%%%%%%%%%%%%%%%%%%%%%%%%%%%%%%%%%%%%%%%%%%%%%%%%%%
%
% MAIN
%
%%%%%%%%%%%%%%%%%%%%%%%%%%%%%%%%%%%%%%%%%%%%%%%%%%%%%%%%%%%%%%%%%%%%%%%%%%%%%%%%%%%%%%%%%%%%%

\section{Main results}
\label{sec:main}
\setcounter{equation}{0}
%%%%%%%%%%%%%%%%%%%%%%%%%%%%%%%%%%%%%%%%%%%%%%%%%%%%%%%%%%%%%%%%%%%%%%%%%%%%%%%%%%%%%%%%%%%%%

\subsection{Preliminaries}
\label{subsec:prel}

We assume that $\Omega$ is a smooth bounded domain of $\RR^3$ with smooth boundary $\Gamma:=\de\Omega$ and we denote $|\Omega|$ for its Lebesgue measure. Let $X$ be a (real) Banach or Hilbert space, whose norm is denoted
by $\|\cdot\|_X$. $X'$ indicates the
dual space of $X$ and $\langle \cdot,\cdot\rangle$ denotes the
corresponding duality product. The boldface letter $\mathbf{X}$
stands for the
vectorial space $X^3$ endowed with the product structure. 
We denote by $L^p(\Omega)$ and $W^{k,p}(\Omega)$, $k\in \mathbb{N}$ and
$p\in [1,+\infty]$, the usual Lebesgue spaces and Sobolev spaces of
real measurable functions on the domain $\Omega$.  We
indicate by $H^k(\Omega)$ the Hilbert spaces
$W^{k,2}(\Omega)$ with respect to the scalar product
$\langle u,v \rangle_{k}=\sum_{|\alpha|\leq
k}\int_{\Omega} D^\alpha u(x) D^\alpha v(x) \, \dix$ ($\alpha=(\alpha_1,\alpha_2,\alpha_3)$ being a
multi-index) and the induced norm $\|u\|_{H^k(\Omega)}=\sqrt{\langle u,u \rangle_k}$.
Given an interval $J$ of $\mathbb{R}^+$, we introduce
the function space $L^p(J;X)$ with $p\in [1,+\infty]$, which
consists of Bochner measurable $p$-integrable
functions with values in the Banach space $X$.

Set $$H := L^2(\Omega),\quad V := H^1(\Omega).$$ 
Then we denote
by $(\cdot,\cdot)$ the standard scalar product of~$H$ and by
$\| \cdot \|$ the associated Hilbert norm.
Identifying $H$ with its dual space $H'$ by means
of the above scalar product, we obtain the chain of continuous
and dense embeddings $V\subset H \subset V'$. Let $\bn$ be the exterior
unit normal vector to $\Gamma$. We set
\begin{equation}\label{defiW}
   W := \big\{ u\in H^2(\Omega):~\dn u=0~\text{on $\Gamma$}\big\},
\end{equation}
 such that $W$ is a closed subspace of $H^2(\Omega)$ (and in particular
it inherits its norm). For every $g\in V'$, we denote by $\overline{g}$ the generalized mean value of function $g$
over $\Omega$ such that
\begin{align}
\overline{g}=\frac{1}{|\Omega|}\langle g,1\rangle. \label{defi:oo}
\end{align}
If $g\in L^1(\Omega)$, then $\overline{g}=|\Omega|^{-1}\int_\Omega g \dix$. 
In this paper we will use the Poincar\'{e}--Wirtinger inequality
\begin{equation}
\label{poincare}
\|g-\overline{g}\|\leq C_P\|\nabla g\|,\quad \forall\,
g\in V,
\end{equation}
where $C_P$ is a constant depending only on $n$ and $\Omega$.
We introduce the linear spaces
$$
V_0=\{ u \in V:\ \overline{u}=0\},\quad
L^2_0(\Omega)=\{ u \in H:\ \overline{u}=0\}, \quad
V_0'= \{ u \in V':\ \overline{u}=0 \},
$$
and we consider the realization of the Laplace operator with homogeneous Neumann boundary conditions 
$A\in \mathcal{L}(V,V')$ defined by 
\begin{equation}\label{defiA}
   \duav{Au,v} := \io \nabla u\cdot \nabla v \, \dix,
     \quext{for }\,u,v\in V.
\end{equation}
The restriction of $A$ from $V_0$ onto $V_0'$
is an isomorphism. In particular, $A$ is positively defined on $V_0$
and self-adjoint. We denote its inverse map by $\mathcal{N} =A^{-1}: V_0'
\to V_0$. Note that for every $g\in V_0'$, $u= \mathcal{N} g \in V_0$ is the unique
(in $V_0$) weak solution of the Neumann problem
$$
\begin{cases}
-\Delta u=g, \quad \text{in} \ \Omega,\\
\dn u=0, \quad \ \  \text{on}\ \partial \Omega.
\end{cases}
$$
Besides, we have
\begin{align}
&\langle A u, \mathcal{N} g\rangle = \langle  g,u\rangle, \quad \forall\, u\in V, \ \forall\, g\in V_0',\label{propN1}\\
&\langle  g, \mathcal{N} h\rangle = \langle h, \mathcal{N} g \rangle = \int_{\Omega} \nabla(\mathcal{N} g)
\cdot \nabla (\mathcal{N} h) \, \dix, \quad \forall \, g,h \in V_0'.\label{propN2}
\end{align}
For any $g\in V_0'$, we set $\|g\|_{V_0'}=\|\nabla \mathcal{N} g\|$.
It is well-known that $g \to \|g\|_{V_0'}$ and $
g \to(\|g-\overline{g}\|_{V_0'}^2+|\overline{g}|^2)^\frac12$ are
equivalent norms on $V_0'$ and $V'$,
respectively. Besides, according to Poincar\'{e}'s inequality \eqref{poincare}, we
have that  $g\to (\|\nabla g\|^2+|\overline{g}|^2)^\frac12$ is an
equivalent norm on $V$.
We also report the following standard Hilbert interpolation inequality and elliptic estimates for the Neumann problem:
\begin{align}
\|g\| &\leq \|g\|_{V_0'}^{\frac12} \| \nabla g\|^{\frac12},
\qquad \forall\, g \in V_0,\label{I}\\
\|\nabla \mathcal{N} g\|_{H^{k}(\Omega)}& \leq C \|g\|_{H^{k-1}(\Omega)},
\qquad \forall\, g\in H^{k-1}(\Omega)\cap L^2_0(\Omega),\quad k=1,2.\label{N}
\end{align}

In order to manage the singular terms related to the choice
of the logarithmic potential~\eqref{defiF}, it is convenient
to introduce some additional notation. First of all,
we indicate by $\beta$ the monotone part of $F'$, namely
\begin{equation}\label{defibeta}
  \beta(r) := \frac12 \ln (1 + r) - \frac12 \ln (1 - r),
    \quad r\in (-1,1).
\end{equation}
Hence, according to~\eqref{defif}, it holds
$$f(r)=F'(r)=\beta(r)-\lambda r.$$ 
We also set
\begin{equation}\label{defia}
  a(r) := 2 \beta'(r) = \frac2{1-r^2}, \quad r\in (-1,1),
\end{equation}
so that
$$f'(r)=F''(r)=\frac{a(r)}{2}-\lambda.$$
For further convenience, we also compute
\begin{equation}\label{aprimo}
  a'(r) = 2 \beta''(r) = \frac{4r}{(1-r^2)^2}, \quad
  a''(r) = 2 \beta'''(r) = \frac{4(1+3r^2)}{(1-r^2)^3}, \quad r\in (-1,1).
\end{equation}
It is also worth rewriting \eqref{eq:u}--\eqref{eq:om} as a single
equation. Recalling \eqref{defif} and \eqref{defia},
we have
\begin{equation}\label{eq:uom}
  \mu = \Delta^2 u - \Delta \beta(u)
   - \beta'(u) \Delta u + \beta(u) \beta'(u)
    + (2\lambda - \eta) \Delta u + g(u),
\end{equation}
where 
\begin{equation}\label{defig}
  g(u):= - \lambda u \beta'(u)
   + ( \eta - \lambda) \beta(u)
   + (\lambda^2 - \lambda \eta) u.
\end{equation}
If $\lambda=\eta=0$, then $g\equiv 0$ which is an easier situation. If $\lambda\neq 0$, the function $g=g(r)$ is bounded on any compact set $I\subset (-1,1)$ and diverges as $|r|\nearrow 1$. Furthermore, we note that the function $\beta(r)\beta'(r)$ is monotone and it dominates $g(r)$ near the pure phases $\pm 1$, namely,
\begin{equation}\label{domin}
  \lim_{|r|\nearrow 1} \frac{| \beta(r) \beta'(r) |}{ | g(r) | }
   = +	\infty.
\end{equation}

Noting that $\Delta \beta(u)=\beta'(u)\Delta u+\beta''(u)|\nabla u|^2$, 
one can rewrite equation \eqref{eq:uom} in several
alternative forms, for instance,
\begin{align}\label{eq:uom1}
  \mu & = \Delta^2 u
    - 2 \Delta \beta(u)
    + \beta''(u) | \nabla u |^2
    + \beta(u) \beta'(u)
    + (2\lambda - \eta) \Delta u + g(u),\\
 \label{eq:uom2}
  \mu & = \Delta^2 u
    - 2 \beta'(u) \Delta u
    - \beta''(u) | \nabla u |^2
    + \beta(u) \beta'(u)
    + (2\lambda - \eta) \Delta u + g(u).
\end{align}
It is however necessary to remark that the above expressions are completely
equivalent only as far as {\sl smooth solutions}\/ are considered. Because
we will only deal with {\sl weak solutions}, in that framework the 
equivalence will be lost and it will be necessary to choose
the most appropriate of the above expressions (namely, \eqref{eq:uom1})
and restate it in a variational form. Using \eqref{defia}, \eqref{aprimo}, we also have 
\begin{align}
 \label{eq:uom3}
  \mu & = \Delta^2 u
    - a(u) \Delta u
    - \frac{a'(u)}2 | \nabla u |^2
    + \beta(u) \beta'(u)
    + (2\lambda - \eta) \Delta u + g(u).
\end{align}
From \eqref{eq:uom3} 
we observe that the singular diffusion terms in $\mu$ 
have exactly the same shape as in \cite{SP0,SP}. Indeed, we borrowed some notation (in particular the use of the coefficient $a(\cdot)$) from
those papers for later convenience of the reader.

%%%%%%%%%%%%%%%%%%%%%%%%%%%%%%%%%%%%%%%%%%%%%%%%%%%%%%%%%%%%%%%%%%%%%%%%%%%%%%%%%%%%%%%%%%%%%

\subsection{Weak formulation and main results}
\label{subsec:weak}

To begin with, we specify our assumptions on the initial datum:
\begin{equation}\label{hp:init}
  u_0 \in W\quad \text{and}\quad  \beta(u_0) \in H.
\end{equation}
As will be further discussed below, hypothesis \eqref{hp:init} corresponds
to the finiteness of the initial energy $\calE(u_0)$ (cf.~\eqref{defiE}).
It is worth observing that \eqref{hp:init} implies $-1 < u_0(x) < 1$
a.e.~in $\Omega$, whence we have in particular  
\begin{equation}\label{hp:mean}
  -1 < m < 1, \quext{where }\,
    m:= \overline{u_0}.
\end{equation}
Indeed, the case $m=1$ (and similarly happens for $m=-1$), implying $u_0=1$
a.e., is incompatible with the hypothesis $\beta(u_0) \in H$.

Next, we present our basic concept of weak solutions:
\bede\label{def:weak}
 A couple $(u,\mu)$ is called a weak solution to problem \eqref{eq:mu}--\eqref{ini} over the time interval $(0,T)$ provided that:\\[2mm]
 (A)~~The following regularity conditions are satisfied:
 \begin{align}\label{reg:u}
   & u \in H^1(0,T;V') \cap L^\infty(0,T;W) \cap L^4(0,T;H^3(\Omega)),\\
  \label{reg:bu1}
   & \beta(u) \in L^\infty(0,T;H) \cap L^4(0,T;V),\\
  \label{reg:bu2}
   & \beta(u)\beta'(u) \in L^2(0,T;L^1(\Omega)),\\
  \label{reg:bu3}
   & \beta''(u) |\nabla u|^2 \in L^2(0,T;L^1(\Omega)),\\
  \label{reg:mu}
   & \mu \in L^2(0,T;V).
 \end{align}
 (B)~~The following weak counterparts of equations\/ \eqref{eq:mu} and \eqref{eq:uom1}
 hold a.e.~in~$(0,T)$:
 \begin{align}\label{E:mu}
  & \partial_t u + A \mu = 0 \quext{in }\,V',\\
  \label{E:u}
   & \mu = A^2 u + 2 A \beta(u)
    + \beta''(u) | \nabla u |^2
    + \beta(u) \beta'(u)
    - (2\lambda - \eta) A u + g(u)
     \quext{in }\,V' + L^1(\Omega),
 \end{align}
 where $A$ is defined by \eqref{defiA}.\\[2mm]
 (C)~~The initial condition is satisfied in the following sense:
 \begin{equation}\label{E:init}
   u|_{t=0} = u_0 \quext{a.e.~in }\,\Omega.
 \end{equation}
\edde
\beos\label{rem:reg0}
It is worth observing that relation~\eqref{E:u} is asked to hold
in the space $V' + L^1(\Omega)$
as a natural consequence of the regularity conditions \eqref{reg:u}--\eqref{reg:mu}.
Equivalently, one can rephrase \eqref{E:u} as the variational equality
\begin{align}\no
  ( \mu, \fhi) & = - ( \nabla \Delta u, \nabla \fhi )
    + 2 ( \nabla \beta(u), \nabla \fhi )
    + \io \beta''(u) | \nabla u |^2 \fhi \dix \\
 \label{eq:uom:w}
  & \mbox{}~~~~~
    + \io \big(\beta(u) \beta'(u)+ g(u)\big) \fhi \dix
    + (2\lambda - \eta)(\Delta u, \fhi ),
\end{align}
for almost all $t\in(0,T)$ and any test function $\fhi\in V \cap L^\infty(\Omega)$ (note that $L^1(\Omega) \subset (L^\infty(\Omega))'$). 
Recalling \eqref{eq:om}, we infer from \eqref{reg:u}--\eqref{reg:bu1} that 
$$ 
\omega \in L^\infty(0,T; H)\cap L^4(0,T; V)\quad \text{and}\quad 
 \omega=-\Delta u + f(u)\quext{a.e.~in }\,(0,T)\times\Omega.
$$
Besides, \eqref{reg:u}--\eqref{reg:bu1} imply that $u\in C_w([0,T];W)$, 
$$ -1<u(x,t)<1 \quext{a.e.~in }\,(0,T)\times\Omega, \quad\text{and}\quad 
\|u(t)\|_{L^\infty(\Omega)}\leq 1\quext{for a.a. }\,t\in(0,T).
$$
\eddos

\smallskip

Now we state the main results of this paper. 
\bete\label{thm:main} {\rm (Well-posedness).}
 Let $F$ be determined by\/ \eqref{defiF} and $\lambda,\,\eta\in \RR$ be given. For any initial datum $u_0$ that satisfies \eqref{hp:init}--\eqref{hp:mean}, there exists a unique weak solution $(u,\mu)$ to problem \eqref{eq:mu}--\eqref{ini} in the sense of Definition~\ref{def:weak}. 
 Moreover, let $u\zzu,u\zzd$ be a couple of initial data
 both satisfying\/ \eqref{hp:init} and such that $\overline{u\zzu}=\overline{u\zzd}\in (-1,1)$.
 Then, denoting $u_1$, $u_2$ the
 corresponding pair of weak solutions to problem \eqref{eq:mu}--\eqref{ini},
 the following continuous dependence estimate holds:
 \begin{equation} \label{cont:dep}
   \| u_1(t) - u_2(t) \|_{V'_0}^2 + \itt \| A u_1(s) - A u_2(s) \|^2 \,\dis
    \le \| u\zzu - u\zzd \|_{V'_0}^2 e^{Ct}, \quext{for all }\,t\in [0,T],
 \end{equation}
 where the constant $C>0$ depends only on the assigned parameters of the problem and,
 in particular, is independent of $u\zzu$, $u\zzd$, and of time.
\ente
Next, we can obtain further properties of weak solutions, like parabolic regularization for strictly positive times and
a suitable form of the energy dissipation principle:
\bete\label{thm:furth}{\rm (Regularity and energy identity).}
 Assume that the hypotheses of Theorem~\ref{thm:main} are satisfied. Let $(u,\mu)$ be the corresponding weak solution defined over the generic interval $[0,T]$. 
 Then, for any $\tau \in (0,T)$
 the following additional regularity properties hold:
 \begin{align}\label{reg:u+}
   & u \in W^{1,\infty}(\tau,T;V') \cap H^1(\tau,T;W) \cap L^\infty(\tau,T;H^3(\Omega)),\\
  \label{reg:bu1+}
   & \beta(u) \in H^1(\tau,T;H) \cap L^\infty(\tau,T;V),\\
  \label{reg:bu2+}
   & \beta(u)\beta'(u) \in L^\infty(\tau,T;L^1(\Omega)),\\
  \label{reg:bu3+}
   & \beta''(u) |\nabla u|^2 \in L^\infty(\tau,T;L^1(\Omega)),\\
  \label{reg:mu+}
   & \mu \in L^\infty(\tau,T;V).
 \end{align}
 Moreover, for any $t_1,t_2$ with $0\leq t_1 < t_2 \le T$ we have the
 energy equality:
 \begin{equation}\label{diss:E2}
   \calE(u(t_2)) + \int_{t_1}^{t_2} \| \nabla \mu(s) \|^2 \, \dis = \calE(u(t_1)).
 \end{equation}
\ente
%
%\beos\label{oss:ener}
%\III 
%If we allow $t_1=0$ in \eqref{diss:E2}, then the equality should be replaced by an energy inequality,
%
% \begin{equation}\label{diss:E2b}
%   \calE(u(t)) + \int_{0}^{t} \| \nabla \mu(s) \|^2 \, \dis \leq \calE(u_0),\quad \text{for all}\ t\in (0,T].
% \end{equation}
%\EEE
%\eddos
%

Finally, we characterize the long-time behavior
of solutions. In view of the global well-posedness result (Theorem~\ref{thm:main})
and the smoothing property (Theorem~\ref{thm:furth}) of weak solutions, it is worth
expecting the existence of a global attractor.
The statement of this natural property requires, however, the introduction of
some further machinery, especially related to the characterization of the
most convenient phase space for the dynamical process associated to the evolution system.
Indeed, looking at the regularity \eqref{hp:init} imposed on the initial datum, we have to take
the singular function~$\beta(u)$ into account. Moreover, we need to consider the mass conservation constraint~\eqref{st:mean}. This leads to the following
\bede\label{def:phase}
 For any given $m\in (-1,1)$, we set
 \begin{equation}\label{defiXm}
   \calX_m:=\big\{ v\in W:~\beta(v) \in H,~\overline{v} = m\big\}.
 \end{equation}
The distance on the phase space $\calX_m$ is defined as follows 
 \begin{equation}\label{distXm}
   \dist_{\calX}(v_1,v_2):= \| v_1 - v_2 \|_W + \| \beta(v_1) - \beta(v_2) \|.
 \end{equation}
\edde
\noindent Clearly, $ \calX_m$ cannot have a linear structure. Nevertheless,
following the lines, e.g., of \cite[Lemma~3.8]{RS}, one can easily
show that $\calX_m$ is a \textit{complete metric space}. Then we have the following result:
\bete\label{thm:attra}{\rm (Global attractor).}
 Assume that the hypotheses of Theorem~\ref{thm:main} are satisfied. Let in
 particular $\overline{u_0}=m$ with assigned $m\in(-1,1)$. Then, the global weak
 solutions to problem \eqref{eq:mu}--\eqref{ini} 
 generate a dynamical process $S(t)$ on the phase
 space~$\calX_m$ which admits a compact global attractor denoted by $\calA_m$.
 Moreover,  we have
 \begin{equation}\label{reg:attra}
   \| u \|_{H^3(\Omega)}
    + \| \beta(u) \|_{V}
    + \| \beta(u)\beta'(u) \|_{L^1(\Omega)}
   \le C_m, \quad \perogni u\in \calA_m.
 \end{equation}
 where the constant $C_m>0$ depends on $\Omega$, $\eta$, $\lambda$, and on the initial datum $u_0$ only through the conserved quantity $m$.
\ente
\beos\label{rem:reg}
 The regularity \eqref{reg:attra} that we can prove for the elements
 of the attractor $\calA_m$ turns out to coincide, as expected, with that
 provided by the parabolic smoothing estimates
 (cf.~\eqref{reg:u+}--\eqref{reg:mu+}). It is then worth discussing
 whether this regularity is optimal or additional properties
 could be proved. 
 
 Generally speaking, when singular potentials are involved,
 the basic regularity threshold for Cahn--Hilliard-like systems
 is linked to the validity of the so-called {\sl strict separation property}:
 \begin{equation}\label{separ}
   \esiste \delta\in (0,1):~~
    -1 + \delta \le u(t,x) \le 1 - \delta
    \quext{for a.e.~}\,x\in \Omega.
 \end{equation}
 Whenever \eqref{separ} holds at some time $t$, then $u(t)$ stays uniformly away in $\Omega$ from the ``singular values'' $\pm 1$. As a consequence,
 singular terms can be treated as smooth functions, which will give rise to
 a further gain of regularity. However, establishing \eqref{separ} for
 some, or all, times $t$ is often a nontrivial question.
 For instance, to the best of our knowledge, \eqref{separ} for any positive time is an open issue for the standard fourth-order Cahn--Hilliard
 system with logarithmic potential, at least in three dimensions (see \cite{MZ}, also \cite{GGM}). On the other hand, \eqref{separ} for any $t>0$ can be proved for
 the Cahn--Hilliard equation with singular diffusion studied in \cite{SP}, which is closely related to our model.
 Unfortunately, the proof of \cite{SP} cannot be reproduced here
 since the argument therein is based on a rather sharp use of \textit{a priori}
 estimates of ``entropy'' type, which do not seem to be applicable
 to sixth order problems (actually the equation addressed in \cite{SP}
 is of fourth order in space). Hence, since the available smoothing estimates \eqref{reg:u+}--\eqref{reg:mu+} appear too weak to guarantee
 \eqref{separ} for the three dimensional case, the question whether additional
 regularity could be proved for our model remains as an open (and likely nontrivial) issue. 
\eddos

Nevertheless, when the spatial dimention is lower than three, we are able to prove the instantaneous strict separation property for weak solutions to problem \eqref{eq:mu}--\eqref{ini}: 
\begin{prop}\label{prop:sepa} {\rm (Separation from pure states $\pm 1$).}
Assume that $\Omega$ is a smooth bounded domain in~$\RR^d$ ($d=1,2$) and all the  other hypotheses of Theorem~\ref{thm:main} are satisfied. Problem \eqref{eq:mu}--\eqref{ini} admits a unique global weak solution $(u,\mu)$ on the generic interval $[0,T]$. Moreover, for any $\tau \in (0,T)$, there exists $\delta\in (0,1)$ such that 
\begin{equation}\label{separ1}
   \|u(t)\|_{L^\infty(\Omega)} \le 1 - \delta,\quad \forall\, t\in [\tau,\,T].
 \end{equation}
\end{prop}

%%%%%%%%%%%%%%%%%%%%%%%%%%%%%%%%%%%%%%%%%%%%%%%%%%%%%%%%%%%%%%%%%%%%%%%%%%%%%%%%%%%%%%%%%%%%%
%
% A PRIORI
%
%%%%%%%%%%%%%%%%%%%%%%%%%%%%%%%%%%%%%%%%%%%%%%%%%%%%%%%%%%%%%%%%%%%%%%%%%%%%%%%%%%%%%%%%%%%%%

\section{\textit{A priori} estimates}
\label{sec:apriori}
\setcounter{equation}{0}

In this section, we will derive a number of \textit{a priori} estimates
for the solution to our problem \eqref{eq:mu}--\eqref{ini}. The estimates will be performed in
a formal way working directly on the ``original'' form of the
system without referring to any explicit approximation
or regularization scheme. In particular, we will always
assume sufficient regularity in order for our computations to make
sense. In such a setting, the various reformulations of \eqref{eq:uom}
(namely, \eqref{eq:uom1}, \eqref{eq:uom2}, \eqref{eq:uom3})
may be assumed to be equivalent to each other. In the next section,
we will discuss a possible approximation scheme (we cannot speak
of ``regularization'' because in fact we will add a further singular
term in the system). Moreover, we will see that such a scheme
is compatible with the estimates derived below up to minor modifications
in the notation and to the management of some additional terms.

In what follows, $\kappa> 0$ and $C\ge 0$ are suitable constants
whose value may vary on occurrence. The
values of $\kappa,\,C$ will be allowed to depend only on the data
of the problem. So, no dependence is allowed on any hypothetical
approximation or regularization parameter. At this stage, dependence
on the ``final time'' $T>0$ is admitted; time-uniform estimates
will be discussed in Section \ref{sec:long}. 
Note, finally, that $\kappa$ is asked to
be strictly positive in view of the fact that it will appear in
estimates from below. 

\medskip

\noindent%
{\bf Energy estimate.}~~%
The first \textit{a priori} information on the solution can be obtained
by reproducing the variational principle. Starting from the equations, this corresponds to
testing \eqref{eq:mu} by $\mu$, \eqref{eq:u} by $\partial_t u$,
integrating over $\Omega$, and taking the difference of the
resultants. Then, using integration by parts and the boundary
condition \eqref{bd}, it follows from \eqref{eq:om} that
\begin{equation}\label{conto:11}
  \io \big( \omega f'(u) \partial_t u  - \Delta \omega \partial_t u) \dix
   = \io \big( \omega f'(u) \partial_t u  - \omega \Delta \partial_t u) \dix
   = \frac12 \ddt \| \omega \|^2.
\end{equation}
As a consequence, we can easily recover the {\sl basic energy law}:
\begin{equation}\label{st:energy}
  \ddt \calE(u) + \| \nabla \mu \|^2
   = 0, \quad \forall\, t\geq 0,
\end{equation}
which implies that the energy functional $\mathcal{E}(u(t))$ is non-increasing in time, i.e., $\mathcal{E}(u(t))\leq \mathcal{E}(u_0)<+\infty$ due to our assumption on the initial datum.

Next, we show that $\mathcal{E}(u)$ is bounded from below and actually
enjoys the $H^2$-coercivity.
Noting that the function $F(r)$ is bounded on the interval $[-1,1]$
and that $f'(r)\ge -\lambda$ for all $r\in(-1,1)$,
then using integration by parts and the boundary
condition \eqref{bd}, we have
\begin{align}
 \no \calE(u)
 &= \frac12 \| \Delta u \|^2
   + \frac12 \| f(u) \|^2-\int_\Omega \Delta u f(u) \dix +\eta\int_\Omega \left(\frac12|\nabla u|^2+F(u)\right)\dix\\
 \no
 &= \frac12 \| \Delta u \|^2
   + \frac12 \| f(u) \|^2+\int_\Omega f'(u)|\nabla u|^2 \dix +\eta\int_\Omega \left(\frac12|\nabla u|^2+F(u)\right)\dix\\
 \no
  & \ge \frac12 \| \Delta u \|^2
   + \frac12 \| f(u) \|^2 +\left(\frac{\eta}{2}-\lambda\right)\|\nabla u\|^2
   +\eta\int_\Omega F(u)\dix \\
 \no
  & \ge \frac12 \| \Delta u \|^2
   + \frac12 \| f(u) \|^2 - C \big( 1 + \| \nabla u \|^2 \big)\\
 \no
 & \ge \frac12 \| \Delta u \|^2
   + \frac12 \| f(u) \|^2
   - C \big( 1 + \| \Delta u \| \| u \| \big)\\
  \no
 & \ge \frac14 \| \Delta u \|^2
   + \frac12 \| f(u) \|^2
   - C(1+\|u\|^2)\\
   \label{coerc:en}
 & \ge \frac14 \| \Delta u \|^2
   + \frac12 \| f(u) \|^2
   - C(\eta, \lambda, |\Omega|).
\end{align}
In the above estimate, we have used the fact that, as long as
$u$ is a function in $H^2(\Omega)$ with finite energy $\mathcal{E}(u)$, then it is necessarily
$-1 < u < 1$ almost everywhere in $\Omega$, otherwise $f(u)$ could not
lie in $L^2(\Omega)$.

We recall that $u$ satisfies the mass-conservation property \eqref{st:mean}.
Combining it with \eqref{st:energy}, using \eqref{coerc:en}
and the Poincar\'e-Wirtinger inequality, we then obtain a control
of the full $W$-norm of $u$. In particular, we have the
following \textit{a priori} estimates:
\begin{align}\label{st:11}
 & u \in H^1(0,T;V') \cap L^\infty(0,T;W),\\
  \label{st:12}
 & \nabla \mu \in L^2(0,T;H),\\
  \label{st:13}
 & \beta(u) \in L^\infty(0,T;H).
\end{align}
Note in particular that the first bound in \eqref{st:11} follows from the control
of $\nabla \mu$ \eqref{st:12} together with the continuity properties of the
operator $A:V\to V'$ (cf.~\eqref{defiA}).

\medskip

\noindent%
{\bf Second estimate.}~~%
 Let us test \eqref{eq:uom1} by $u-\overline{u}$ to obtain
\begin{align}\no
  & \io \mu( u - \overline{u} ) \dix
   = \| \Delta u \|^2
   + \io \big( 2 \beta'(u) + ( u - \overline{u}) \beta''(u) \big) | \nabla u |^2 \dix \\
 \label{co:20}
  & \mbox{}~~~~~ + \io \big( \beta(u)\beta'(u) + g(u) \big)(u-\overline{u}) \dix
   - ( 2\lambda - \eta ) \| \nabla u \|^2.
\end{align}
Let $m\in (-1,1)$ be given. Recalling the monotonicity of the function
$\beta(r)\beta'(r)$ on $(-1,1)$ and the fact $\beta(0)\beta'(0)=0$, 
we have 
\begin{align}
\beta(r)\beta'(r)(r-m)
&= \beta(r)\beta'(r)\Big(r-\frac{|m|+1}{2}\Big) + \beta(r)\beta'(r)\Big(\frac{|m|+1}{2}-m\Big)\nonumber\\
&\geq |\beta(r)\beta'(r)|\Big(\frac{|m|+1}{2}-m\Big) ,\quad \text{for }r\in \Big[\frac{|m|+1}{2},1\Big), \nonumber
\end{align}
\begin{align}
\beta(r)\beta'(r)(r-m)
&= \beta(r)\beta'(r)\Big(r+\frac{|m|+1}{2}\Big) - \beta(r)\beta'(r)\Big(\frac{|m|+1}{2}+m\Big)\nonumber\\
&\geq |\beta(r)\beta'(r)|\Big(\frac{|m|+1}{2}+m\Big) ,\quad \text{for }r\in \Big(-1,\,-\frac{|m|+1}{2}\Big], \nonumber
\end{align}
and
\begin{align}
|\beta(r)\beta'(r)(r-m)|\leq C_m,\quad\text{for }r\in \Big(-\frac{|m|+1}{2},\,\frac{|m|+1}{2}\Big),\nonumber
\end{align}
where the constant $C_m$ is determined by $\frac{|m|+1}{2}$. As a consequence, 
it holds
\begin{align}
\beta(r)\beta'(r)(r-m)\geq \min\left\{\frac{|m|+1}{2}+m,\,\frac{|m|+1}{2}-m\right\} |\beta(r)\beta'(r)| -C_m,\quad \text{for } r\in (-1,1).
\end{align}
Then, using the mass conservation property~\eqref{st:mean}, we can take $m=\overline{u_0}$ and then find constants
$\kappa>0$, $C\ge 0$ depending on $m$ such that
\begin{equation}\label{co:21a}
  \io \beta(u)\beta'(u) (u-\overline{u}) \dix
   \ge \kappa \| \beta(u) \beta'(u) \|_{L^1(\Omega)} - C.
\end{equation}
On the other hand, in view of \eqref{domin}, there exists a $r^*\in(0,1)$ such that $|g(r)|\leq \frac{\kappa}{4}|\beta(r)\beta'(r)|$ for $|r|\in [r^*,1)$, where $\kappa$ 
is the one as in \eqref{co:21a}. Then we have
\begin{align}
&\left|\io g(u)(u-\overline{u}) \dix \right|\nonumber\\
&\quad \leq \left|\int_{\{\Omega:\ |u(x)|\geq r^*\}} g(u)(u-\overline{u}) \dix \right| +  \left|\int_{\{\Omega:\ |u(x)|< r^*\}} g(u)(u-\overline{u}) \dix \right|\nonumber\\
&\quad \leq \frac{\kappa}{2} \| \beta(u) \beta'(u) \|_{L^1(\Omega)}+ C,\nonumber
\end{align}
where $C$ may depend on $r^*$ and $\Omega$. Therefore, we can conclude that 
\begin{equation}\label{co:21}
  \io \big( \beta(u)\beta'(u) + g(u) \big)(u-\overline{u}) \dix
   \ge \frac{\kappa}{2} \| \beta(u) \beta'(u) \|_{L^1(\Omega)} - C.
\end{equation}
Hence, from \eqref{co:20} and \eqref{co:21} we obtain
\begin{align}\no
  & \frac{\kappa}{2} \| \beta(u) \beta'(u) \|_{L^1(\Omega)}
   + \| \Delta u \|^2
   + \io \big( 2 \beta'(u) + ( u - \overline{u}) \beta''(u) \big) | \nabla u |^2 \dix \\
 \no
  & \mbox{}~~~~~
   \le C + (2\lambda - \eta) \| \nabla u \|^2
   + \io \mu ( u - \overline{u}) \dix\\
  \no
  & \mbox{}~~~~~
   \le C +\io  \overline{\mu} ( u - \overline{u}) \dix + \io (\mu - \overline{\mu}) ( u - \overline{u}) \dix \\
 \no
  & \mbox{}~~~~~ \le C+ C \| \nabla \mu \| \| \nabla u \|\\
 \label{co:22}
  & \mbox{}~~~~~ \le C \big(1 + \| \nabla \mu \| \big),
\end{align}
where we have used estimate \eqref{st:11} and the
Poincar\'e--Wirtinger inequality. 
Besides, recalling \eqref{defia} and \eqref{aprimo}, we also have
\begin{align}
 \no
   &  \io \big( 2 \beta'(u) + ( u - \overline{u}) \beta''(u) \big) | \nabla u |^2 \dix \\
 \no
   & \mbox{}~~~~~ = \io \frac{2(1-u \overline{u})}{(1-u^2)^2} | \nabla u |^2 \dix \\
 \no
   & \mbox{}~~~~~
   \ge \io \frac{\kappa_1}{(1-u^2)^2} | \nabla u |^2 \dix \\
 \label{co:23}
   & \mbox{}~~~~~ = \kappa_1\io \beta'(u)^2 | \nabla u |^2 \dix
   = \kappa_1\| \nabla \beta(u) \|^2,
\end{align}
where $\kappa_1=2(1-|\overline{u}|)$ depending on
the conserved value $\overline{u}$ and we used in an essential way the facts that $|\overline{u}|<1$ and $-1\le u \le 1$ almost everywhere in $\Omega$.
Therefore, squaring \eqref{co:22}, integrating in time,
using \eqref{co:23}, and recalling \eqref{st:12},
we arrive at
\begin{align}\label{st:21}
  & \beta(u)\beta'(u) \in L^2(0,T;L^1(\Omega)),\\
 \label{st:22}
  & \nabla \beta(u) \in L^4(0,T;H).
\end{align}
Noting that $\beta'$ is exponentially larger than $|\beta|$ as $|r|\nearrow 1$, then combining
\eqref{st:21} and \eqref{st:22}, we obtain
\begin{equation}\label{st:23}
  \beta(u) \in L^4(0,T;V).
\end{equation}

Next, integrating \eqref{eq:uom1} in space, we deduce that
\begin{equation}\label{st:23a}
  \overline{\mu} = \frac1{|\Omega|} \bigg( \io \beta''(u)|\nabla u|^2\,\dix
   + \io \big( \beta(u) \beta'(u) + g(u) \big) \,\dix \bigg).
\end{equation}
Noting that, by \eqref{st:11}, \eqref{defia}, \eqref{aprimo}, 
and H\"{o}lder's inequality,
\begin{equation}\label{st:23b}
  \bigg| \io \beta''(u)|\nabla u|^2\,\dix \bigg|
   \le \big( \|\nabla u\|^2 
   + 2 \| \nabla \beta(u) \|^2 \big)\leq C \big( 1 + \| \nabla \beta(u) \|^2 \big),
\end{equation}
and recalling \eqref{st:21}, \eqref{st:22}, \eqref{domin}, it is not difficult to arrive at
\begin{equation}\label{st:24a}
  \overline{\mu} \in L^2(0,T),
\end{equation}
which, combined with \eqref{st:12} and the
Poincar\'e--Wirtinger inequality, finally gives
\begin{equation}\label{st:24}
  \mu \in L^2(0,T;V).
\end{equation}

\medskip

\noindent%
{\bf Third estimate.}~~%
Testing \eqref{eq:uom} by $-\Delta u$, after integration by parts, we obtain
\begin{align}\no
  & \| \nabla \Delta u \|^2
   + \io \beta'(u) | \Delta u |^2 \dix
   + \io \big( \beta'(u)^2 + \beta(u)\beta''(u) + g'(u) \big)
    | \nabla u |^2 \dix \\
 \no
  & \mbox{}~~~~~ =  \io \nabla \Delta u \cdot\nabla \beta(u) \dix
   + (2\lambda - \eta) \| \Delta u \|^2
   + \io \nabla \mu \cdot \nabla u \dix \\
 \label{co:31}
  & \mbox{}~~~~~ \le \frac12 \| \nabla \Delta u \|^2
   + \frac12 \| \nabla \beta(u) \|^2
   + C \big( 1 + \| \nabla \mu \| \big),
\end{align}
where we also used the estimate \eqref{st:11} and the Cauchy--Schwarz inequality.

Recalling \eqref{defia} and \eqref{aprimo}, we have
\begin{align}
 & \beta(r)\beta''(r)
    = \frac{r}{(1-r^2)^2} \ln \frac{1 + r}{1 - r}, \label{co:32b}
 \\
 & g'(r) = - \lambda r \beta''(r)+ (\eta - 2\lambda) \beta'(r)
      + \lambda^2 - \lambda\eta
  = - \frac{2\lambda r^2}{(1-r^2)^2}
    + \frac{\eta - 2\lambda}{1-r^2}    
    + \lambda^2 - \lambda\eta.\label{gprimo}
\end{align}
Then for $-1<r<1$, it follows that
\begin{align}
   \beta(r)\beta''(r) + g'(r)
   & = \frac{r}{(1-r^2)^2}\ln\frac{1+r}{1-r}
    - \frac{2\lambda}{(1-r^2)^2}
    + \frac{\eta }{1-r^2}
    + \lambda^2 - \lambda\eta \nonumber\\
 &\ge  \frac{r}{(1-r^2)^2}\ln\frac{1+r}{1-r}
   - \frac{2|\lambda|+2|\eta|}{(1-r^2)^2}  - |\lambda\eta|.
    \label{co:32a}
\end{align}
Since $r \ln\frac{1+r}{1-r}$ is an even function on $(-1,1)$ that is strictly increasing on $(0,1)$,
there exists a $r^*\in [0,1)$ such that $ r^* \ln\frac{1+r^*}{1-r^*}= 4|\lambda|+4|\eta|$ and
$r \ln\frac{1+r}{1-r}\geq 4|\lambda|+4|\eta|$ if $|r|\in[r^*,1)$. Hence, we have in particular
\begin{equation}\label{co:33}
  \beta(r)\beta''(r) + g'(r)
  \ge \frac{1}{2}\frac{r}{(1-r^2)^2}\ln\frac{1+r}{1-r}
   - \frac{2|\lambda|+2|\eta|}{(1-(r^*)^2)^2}  - |\lambda\eta|\quad \text{for}\ r\in (-1,\, 1).
  \end{equation}
As a consequence, it follows from \eqref{co:31} and \eqref{co:33} that  
\begin{align}\label{co:31a}
    & \| \nabla \Delta u \|^2
     \le  \| \nabla \beta(u) \|^2 + C \big( 1 + \| \nabla \mu \| \big).
\end{align}
Thus we can take the
square of \eqref{co:31a} and subsequently integrate the result
in time. Recalling \eqref{st:11}, \eqref{st:12} and
\eqref{st:22}, we conclude the estimate
\begin{equation}\label{st:31}
  u \in L^4(0,T;H^3(\Omega)).
\end{equation}

\medskip

\noindent%
{\bf Fourth estimate.}~~%
Multiplying \eqref{eq:uom} by $\beta(u)$ and integrating
over~$\Omega$, we obtain that 
\begin{align}\no
  & \io \big( 2 \beta'(u)^2 + \beta(u)\beta''(u) \big) | \nabla u |^2\dix
   + \io \beta(u) \big( \beta(u) \beta'(u) +  g(u) \big) \dix\\
 \no
  & \mbox{}~~~~~ = \io \nabla \Delta u \cdot \nabla \beta(u) \dix
   + \io \beta(u) \big( \mu - (2\lambda - \eta) \Delta u \big) \dix\\
 \label{co:41}
  & \mbox{}~~~~~ \le \frac12 \| \nabla \Delta u \|^2
   + \frac12 \| \nabla \beta(u) \|^2
   + C \| \beta(u) \| \big( \| \mu \| + \| \Delta u \| \big).
\end{align}
Using \eqref{defibeta}, \eqref{defia} and \eqref{domin} again, by a direct calculation we can check that 
\begin{equation}\label{co:42}
  \beta(r)^2 \beta'(r) + \beta(r) g(r)
   \ge \frac12 \beta(r)^2 \beta'(r) - C
    = \frac{1}{8(1-r^2)}\ln^2\frac{1+r}{1-r} - C,
\end{equation}
for $r\in (-1,1)$, where $C$ is independent of $r$. 
Hence, computing $\beta(u)\beta''(u)$ using~\eqref{co:32b}, we get
\begin{align}
 \no
  & \io \big( 2 \beta'(u)^2 + \beta(u)\beta''(u) \big) | \nabla u |^2\dix
    + \io \big( \beta(u)^2 \beta'(u) + \beta(u) g(u) \big) \dix\\
 \no
  & \mbox{}~~~~~
    \ge  \io \frac{2 }{(1-u^2)^2} | \nabla u |^2 \dix 
    +\io  \frac{u}{(1-u^2)^2}\left(\ln\frac{1+u}{1-u}\right) | \nabla u |^2\dix \\
\label{co:42a}
  & \mbox{}~~~~~\quad
   + \io \frac{1}{8(1-u^2)}\ln^2\frac{1+u}{1-u} \dix -C.
\end{align}
Recalling \eqref{st:11}, \eqref{st:13}, \eqref{st:23}, \eqref{st:24}
and \eqref{st:31},  we infer from \eqref{co:41} that
\begin{equation}\label{co:41b}
  \io  \frac{u}{(1-u^2)^2}\left(\ln\frac{1+u}{1-u}\right)| \nabla u |^2 \dix
   + \io \frac{1}{(1-u^2)}\ln^2\frac{1+u}{1-u} \dix
   \le \Lambda(t),
\end{equation}
where $\Lambda$ is a nonnegative scalar function of time
satisfying
\begin{equation}\label{co:41c}
  \Lambda(t) \in L^2(0,T).
\end{equation}

Below we prove that the above estimate \eqref{co:41b} yields uniform integrability of the nonlinear terms of \eqref{eq:uom1} involving derivatives of $\beta$. To this aim, we
define, for $-1<u<1$,
\begin{align}\label{co:43}
   &A(u,\nabla u) := \beta''(u) | \nabla u |^2
    = \frac{2u}{(1-u^2)^2} | \nabla u |^2,\\
 \label{co:44}
   &B(u)  := \beta(u) \beta'(u) = \frac1{2(1-u^2)}\ln\frac{1+u}{1-u},
\end{align}
and for $r\ge 0$,
\begin{equation}\label{defiM}
  M(r) := r \ln^{\frac12} (1 + r), \quext{so that }\,
   \lim_{r\nearrow +\infty} \frac{M(r)}r = +\infty.
\end{equation}
It is not difficult to verify that
\begin{equation}
 \no
  \lim_{|r|\nearrow 1} \frac{\displaystyle{\ln\left(1+\left|\frac{1}{1-r^2}\ln\frac{1+r}{1-r}\right|\right)}}
  {\displaystyle{\ln^2\frac{1+r}{1-r}}}=0.
\end{equation}
Then it follows from \eqref{co:22}, \eqref{co:41b} and \eqref{co:44} that
\begin{align}
 \no
  \io M(|B(u)|) \dix
   &
   \le C \io \frac{1}{(1-u^2)}\left|\ln\frac{1+u}{1-u} \right| \bigg( 1 + \left|\ln\frac{1+u}{1-u}\right| \bigg)\dix \\
 \label{st:45}
   &   \le C\big(1+\|\nabla \mu\|+\Lambda(t)\big).
\end{align}

Next, we deal with the term depending on the gradient $\nabla u$, which is a little
bit more tricky.
Let us set, still for $r\ge 0$,
\begin{equation}\label{defiN}
  N(r) := r \ln\ln (e^4 + r), \quext{so that }\,
   \lim_{r\nearrow +\infty} \frac{N(r)}r = +\infty.
\end{equation}
We observe that, for any $x,y\ge 0$,
\begin{align}
 \no
 \ln\ln (e^4 + xy)
   & \le \ln \ln (e^2+x)(e^2+y)\\
 \no
   & = \ln \big( \ln (e^2+x) + \ln (e^2+y) \big)\\
 \label{co:log}
   & \le \ln \ln (e^2+x) + \ln \ln (e^2+y),
\end{align}
thanks to elementary properties of the logarithm. In particular,
we used that $\ln(p+q) \le \ln p + \ln q$ for all $p,q\ge 2$.
Let us now estimate
\begin{align} \no
  \io N(|A(u,\nabla u)|) \dix
   & = \io \frac{2|u|}{(1-u^2)^2} | \nabla u |^2
    \ln\ln\bigg( e^4 + \frac{2|u|}{(1-u^2)^2} | \nabla u |^2 \bigg) \dix \\
 \no
   & \le \io \frac{2|u|}{(1-u^2)^2} | \nabla u |^2
    \ln\ln\bigg( e^2 + \frac{2|u|}{(1-u^2)^2} \bigg) \dix \\
 \no
  & \mbox{}~~~~~
   + \io \frac{2|u|}{(1-u^2)^2} | \nabla u |^2
    \ln\ln\big( e^2 + | \nabla u |^2 \big) \dix\\
  \label{co:46}
   &    =: I_1 + I_2.
\end{align}
Then, from the fact
\begin{equation}
 \no
  \lim_{|u|\nearrow 1} \frac{\displaystyle{\ln\ln\left(e^2+\frac{2|u|}{(1-u^2)^2}\right)}}
  {\displaystyle{\left|\ln\frac{1+u}{1-u}\right|}}=0,
\end{equation}
it is immediate to check from \eqref{co:21}, \eqref{co:23} and \eqref{co:41b} that
\begin{align}
  I_1 &\le  C \io | \nabla u |^2 \frac{2|u|}{(1-u^2)^2}
       \left(1+\left|\ln \frac{1+u}{1-u} \right|\right)\dix \nonumber\\
      &\leq C\|\nabla \beta(u)\|^2+ C\io  \frac{u}{(1-u^2)^2}\left(\ln\frac{1+u}{1-u}\right) | \nabla u |^2\dix \nonumber\\
      &\le C \big(1 + \|\nabla \mu\|+\Lambda(t) \big).
      \label{co:47}
\end{align}
On the other hand, to control $I_2$, we recall that if
$\phi:\RR\to \RR\cup\{+\infty\}$ is a convex and lower semicontinuous function
and $\phi^*$ denotes its convex conjugate function, then
for any $R,S\in \RR$ it holds (see, e.g., \cite{Br})
\begin{equation} \label{conj}
  RS \le \phi(R) + \phi^*(S),
\end{equation}
where the \rhs\ may possibly be $+\infty$ for some $r,s$. 
Now we apply the abstract inequality \eqref{conj} with
the following choices:
\begin{align}  \label{conj:1}
  & \phi(R) = (1 + R) \ln (1 + R), \qquad
   \phi^*(S) = e^{S-1} - S \le e^S~~\text{for }\,S>0,\\
 \label{conj:2}
  & R = \frac{2|u|}{(1-u^2)^2}\geq 0, \qquad
    S = \ln\ln\big( e^2 + | \nabla u |^2 \big)>0.
\end{align}
Using the simple fact
\begin{equation}
 \no
  \lim_{|r|\nearrow 1} \frac{\displaystyle{\left|\left(1+\frac{2|r|}{(1-r^2)^2}\right)\ln\left(1+\frac{2|r|}{(1-r^2)^2}\right)\right|}}
  {\displaystyle{\left|\frac{r}{(1-r^2)^2}\ln\frac{1+r}{1-r}\right|}}
  =4,
\end{equation}
we infer from \eqref{conj:1}, \eqref{conj:2} and the estimates \eqref{st:11}, \eqref{co:41b} that
\begin{align}
I_2&\leq \int_\Omega | \nabla u |^2 \left[\left(1+\frac{2|u|}{(1-u^2)^2}\right)\ln\left(1+\frac{2|u|}{(1-u^2)^2}\right)+\ln\big( e^2 + | \nabla u |^2 \big)\right] \dix\nonumber\\
  &\leq \int_\Omega | \nabla u |^2 \left(C + \frac{5u}{(1-u^2)^2}\ln\frac{1+u}{1-u}\right) \dix+ \int_\Omega | \nabla u |^2 \ln\big( e^2 + | \nabla u |^2 \big) \dix\nonumber\\
    &\leq C\big(1+\Lambda(t)\big).
    \label{co:47b}
\end{align}
As a result, we obtain from \eqref{co:47}, \eqref{co:47b} that
\begin{align}
  \io N(|A(u,\nabla u)|) \dix
   & \le C \big(1 + \|\nabla \mu\|+ \Lambda(t) \big).
   \label{co:48}   
\end{align}

Collecting the estimates \eqref{st:45}, \eqref{co:48} and keeping \eqref{st:24}, \eqref{co:41c} in mind, we arrive at the following estimate
\begin{equation}\label{st:41}
  \| M(|B(u)|) \|_{L^2(0,T;L^1(\Omega))}
   + \| N(|A(u,\nabla u)|) \|_{L^2(0,T;L^1(\Omega))}
   \le C.
\end{equation}
This is a crucial estimate that we need for the purpose of proving existence of a weak solution.

%%%%%%%%%%%%%%%%%%%%%%%%%%%%%%%%%%%%%%%%%%%%%%%%%%%%%%%%%%%%%%%%%%%%%%%%%%%%%%%%%%%%%%%%%%%%%
%
% PROOF OF EXISTENCE
%
%%%%%%%%%%%%%%%%%%%%%%%%%%%%%%%%%%%%%%%%%%%%%%%%%%%%%%%%%%%%%%%%%%%%%%%%%%%%%%%%%%%%%%%%%%%%%

\section{Existence of weak solutions}
\label{sec:exi}
\setcounter{equation}{0}

In this section we shall provide the proof of Theorem~\ref{thm:main}
for the part that concerns the existence of weak solutions.

%%%%%%%%%%%%%%%%%%%%%%%%%%%%%%%%%%%%%%%%%%%%%%%%%%%%%%%%%%%%%%%%%%%%%%%%%%%%%%%%%%%%%%%%%%%%%

\subsection{Weak sequential stability}
\label{subsec:stab}

Our strategy of proof is as follows. As a first step, we assume to have a sequence of sufficiently smooth solutions $(u_n,\mu_n)$ to problem~\eqref{eq:mu}--\eqref{ini}
complying with the uniform \textit{a priori} estimates with respect to~$n$ over the whole reference interval $(0,T)$. 
Then we shall see that, at least, there exists a subsequence converging in a
suitable way to a couple of functions $(u,\mu)$ solving problem~\eqref{eq:mu}--\eqref{ini} in the sense of Definition~\ref{def:weak}. This type of
property, usually referred to as \textit{weak sequential stability}, can be viewed  as an abbreviated method for passing to the limit in a suitable approximation or regularization of the original problem. A possible way to explicitly approximate the system and to construct the sequence $(u_n,\mu_n)$ will be presented in the next subsection. 

Now, we proceed to deduce a number of properties following from the
\textit{a priori} estimates. We point out that all convergence relations stated below are always implicitly assumed to hold up to the extraction of (non-relabelled) subsequences.

First of all, by \eqref{st:11} and \eqref{st:31}, we have
\begin{equation}
\label{conv:11}
  u_n \to u \quext{weakly star in }\,L^\infty(0,T;W)
   \quext{and weakly in }\,H^1(0,T;V') \cap L^4(0,T;H^3(\Omega)).
\end{equation}
Using the Aubin--Lions compactness lemma and the Sobolev embedding theorem, we then obtain, for all $\sigma\in (0,1/2)$,
\begin{equation}\label{conv:12}
  u_n \to u \quext{strongly in }\,C([0,T];H^{2-\sigma}(\Omega))
   \quext{and hence, uniformly in } [0,T]\times\barO.
\end{equation}
Next, from \eqref{st:24}, we also have
\begin{equation}\label{conv:13}
  \mu_n \to \mu \quext{weakly in }\,L^2(0,T;V).
\end{equation}
These convergent results are sufficient to pass to the limit in
equation~\eqref{E:mu} for $(u_n,\mu_n)$. 

The more delicate issue stands, clearly,
in dealing with \eqref{E:u}. Its validity at the $n$-level actually
implies that
\begin{align}\no
  ( \mu_n, \fhi) & = - ( \nabla \Delta u_n, \nabla \fhi )
    + 2 ( \nabla \beta(u_n), \nabla \fhi )
    + \io \beta''(u_n) | \nabla u_n |^2 \fhi \dix \\
 \label{eq:uom:wn}
  & \mbox{}~~~~~
    + \io \big(\beta(u_n) \beta'(u_n)+ g(u_n)\big) \fhi \dix
    + \big( (2\lambda - \eta) \Delta u_n, \fhi \big),
\end{align}
for any $\fhi\in V\cap L^\infty(\Omega)$. Then, in order to take the limit $n\nearrow \infty$ in~\eqref{eq:uom:wn}, we need to manage the nonlinear
terms. First, combining the uniform convergence \eqref{conv:12} and
the bounds \eqref{st:13} and \eqref{st:23}, we have
\begin{equation}\label{conv:13B}
  \beta(u_n) \to \beta(u) \quext{weakly star in }\,
   L^\infty(0,T;H) \cap L^4(0,T;V).
\end{equation}
In particular, the identification of the limit follows from standard
monotonicity method (see, e.g., \cite[Proposition~1.1, Chapter II]{Ba}).
Next, we observe that, as a consequence of \eqref{conv:12},
we have in particular
\begin{equation}\label{conv:14}
  \nabla u_n \to \nabla u \quext{a.e.~in }\,(0,T)\times\Omega,
\end{equation}
whence, also,
\begin{equation}\label{conv:15}
  \beta''(u_n) | \nabla u_n |^2 \to \beta''(u) | \nabla u |^2
    \quext{a.e.~in }\,(0,T)\times\Omega.
\end{equation}
Now, let us discuss the consequences of estimate \eqref{st:41}
applied to the approximating sequence $\{u_n\}$. Indeed,
in view of the fact that the functions $M(|\cdot|)$ and
$N(|\cdot|)$ are convex (as a direct check shows) and coercive
at infinity (as indicated in \eqref{defiM} and \eqref{defiN}),
we may apply the de la Vall\'ee--Poussin criterion
(cf., e.g., \cite[Chapter~2]{DM}) to conclude that
the families $\{A(u_n,\nabla u_n)\}$ and $\{B(u_n)\}$ are
equi-integrable. This fact, combined with the pointwise convergence
\eqref{conv:15}, implies, by Vitali's convergence theorem together with \eqref{st:24} and \eqref{co:41c}, that
\begin{equation}\label{conv:16}
  \beta''(u_n) | \nabla u_n |^2 \to \beta''(u) | \nabla u |^2
    \quext{strongly~in }\,L^2(0,T;L^1(\Omega)).
\end{equation}
Hence, in view of the choice $\fhi \in V \cap L^\infty(\Omega)$,
\begin{equation}\label{conv:17}
  \io \beta''(u_n) | \nabla u_n |^2 \fhi \,\dix \to \io \beta''(u) | \nabla u |^2 \fhi \,\dix
   \quext{a.e.~in }\,(0,T).
\end{equation}
Analogously, Vitali's theorem guarantees that
\begin{equation}\label{conv:16b}
  \beta(u_n) \beta'(u_n) \to \beta(u) \beta'(u)
    \quext{strongly~in }\,L^2(0,T;L^1(\Omega)),
\end{equation}
whence
\begin{equation}\label{conv:17b}
  \io \beta(u_n) \beta'(u_n) \fhi \,\dix \to \io \beta(u) \beta'(u) \fhi \,\dix
   \quext{a.e.~in }\,(0,T).
\end{equation}
The nonlinear term $g(u_n)$ in \eqref{eq:uom:wn}
can be treated in the same way since it is dominated by $\beta(u_n)\beta'(u_n)$
(cf.~\eqref{domin}). As a conclusion, we may take the limit $n\nearrow \infty$ in \eqref{eq:uom:wn}
and recover \eqref{eq:uom:w}. As noted before, this may be equivalently reformulated as
\eqref{E:u}. 

To complete the proof, it is just worth observing that the regularity properties
\eqref{reg:u}--\eqref{reg:mu} are a direct consequence of convergence relations
\eqref{conv:11}, \eqref{conv:13}, \eqref{conv:13B}, \eqref{conv:16} and \eqref{conv:16b}.
Finally, we can pass to the limit in the initial condition and get
back~\eqref{E:init} in view, e.g., of the time-uniform convergence \eqref{conv:12}.

This completes the proof of the existence part of Theorem~\ref{thm:main}.

%%%%%%%%%%%%%%%%%%%%%%%%%%%%%%%%%%%%%%%%%%%%%%%%%%%%%%%%%%%%%%%%%%%%%%%%%%%%%%%%%%%%%%%%%%%%%

\subsection{Approximation scheme}
\label{subsec:appro}

Let us now present a possible approximation of system \eqref{eq:mu}--\eqref{eq:om}.
For any integer $n\geq 3$, we consider
$I_{[-1+1/n,1-1/n]}$, the indicator function of the interval
$[-1+1/n,1-1/n]$ and we associate it to the convex functional
\begin{equation}\label{defId}
  \mathcal{I}_n:V \to [0,+\infty],   \qquad
   \mathcal{I}_n(v):= \io I_{[-1+1/n,1-1/n]}(v(x))\dix.
\end{equation}
In other words, $\mathcal{I}_n(v)$ equals to $0$ if $v\in V$ satisfies
$-1+1/n\le v(x)\le 1-1/n$ a.e.~in~$\Omega$ and is $+\infty$
otherwise. Referring to \cite[Chapter~3]{At} for the definition and basic properties
of Mosco- and graph- convergence, we can prove the following
simple property:
\bele\label{l:mosco}
 The functionals $\mathcal{I}_n$ converge to $\mathcal{I}$ in the sense of Mosco in
 the space $V$, where $\mathcal{I}$ is defined as follows:
 \begin{equation}\label{defI}
   \mathcal{I}:V \to [0,+\infty],   \qquad
    \mathcal{I}(v):= \io I_{[-1,1]}(v(x))\dix.
 \end{equation}
\enle
\begin{proof}
First, we need to prove that, if $v_n$ tends to $v$ weakly in $V$, then
\begin{equation}\label{mosco:1}
  \mathcal{I}(v) \le \liminf_{n\nearrow +\infty} \mathcal{I}_n(v_n).
\end{equation}
To show this, it is sufficient to consider the case when the \rhs\ of \eqref{mosco:1} is
finite. This is indeed equivalent to saying that it is $0$.
Then there exists a (nonrelabelled) subsequence
such that $-1+1/n \le v_n \le 1-1/n$ almost everywhere.
By weak compactness, we can assume
that $v_n$ tends to $v$ strongly in $L^2$ and
a.e.~in $\Omega$, it also follows that $-1\le v \le 1$
almost everywhere. Thus, the \lhs\ of~\eqref{mosco:1}
is $0$, as desired. 

Next, we need to prove that, for any $v\in V$, there
exists a family $\{v_n\}\subset V$ such that $v_n$ tends
to $v$ strongly in $V$ and
\begin{equation}\label{mosco:2}
  \mathcal{I}(v) = \lim_{n\nearrow +\infty} \mathcal{I}_n(v_n).
\end{equation}
Indeed, if $\mathcal{I}(v)=+\infty$, then $|v|>1$ in a set of strictly
positive measure and one can simply take $v_n\equiv v$; otherwise
one may truncate $v$ at the levels $-1+1/n$ and $1-1/n$.
Then, as proved in Lemma~\ref{lem:trunc} below,
$v_n$ tends to $v$ strongly in $V$,
whence the desired property.
\end{proof}

Let $\Gamma_n$ denote the subdifferential
of $\mathcal{I}_n$ in the duality between $V$ and $V'$. Namely,
for any $v\in V$ and $\zeta\in V'$ we set
\begin{equation}\label{defGd}
  \zeta\in \Gamma_n(v) \quext{if and only if}~~
   \duav{\zeta,z-v}_{V',V} + \mathcal{I}_n(v)\le \mathcal{I}_n(z)\quad \perogni z \in V.
\end{equation}
Then by definition, $\Gamma_n$ is a (multivalued) maximal
monotone operator from $V$ to $2^{V'}$. The properties of $\Gamma_n$
have been first described in \cite{brezisart} in the related
case $V=H^1_0(\Omega)$ (corresponding to the homogeneous Dirichlet boundary condition)
and then further characterized in several papers 
(see, e.g., \cite{BCGG,SP0}). Here we just recall that,
if $v\in V$ and $\zeta\in \Gamma_n(v)$, then
$\zeta$, beyond lying (by definition) in $V'$, can also be
interpreted as a measure. Namely, there exists a
Borel measure $\nu$ on $\barO$ such that (cf. \cite[Proposition 2.1]{SP0})
\begin{equation}\label{interpr}
  \duav{\zeta,\fhi}_{V',V} = \ibaro \fhi \deriv\!\nu
    \quad \perogni \fhi\in V \cap C(\barO).
\end{equation}
Note that, once the domain $\Omega$ is smooth, the space $V \cap C(\barO)$
is dense both in $V$ and in $C(\barO)$ (hence the measure $\nu$ is
univocally defined). It may actually happen that $\nu$ is partially
supported on the boundary (and consequently it is essential to integrate over
$\barO$). Some further related properties will be recalled in the estimates
below.

With these preliminaries at hand, we now introduce our approximate system depending on the parameter $n$
(probably it is not proper to speak of a  ``regularized problem'', in the
sense that the system below contains in fact an additional singular term):
\begin{align}\label{eq:mu:d}
  & \partial_t u_{n} + A \mu_n = 0
   \quext{in }\,V',\\
 \label{eq:u:d}
  & \mu_n = A^2 u_n
    + \zeta_n
    + 2 A \beta(u_n)
    + \beta''(u_n) | \nabla u_n |^2
    + \beta(u_n) \beta'(u_n)
    - (2\lambda - \zeta) A u_n + g(u_n) \quext{in }\,V',\\
 \label{eq:eta:d}
  & \zeta_n \in \Gamma_n(u_n) \quext{in }\,V'.
\end{align}
It turns out that dealing
with the system \eqref{eq:mu:d}--\eqref{eq:eta:d}
is indeed simpler than handling the original system \eqref{eq:mu}--\eqref{eq:om}.
The key point is that, once a pair $(u_n,\mu_n)$ solves the system \eqref{eq:mu:d}--\eqref{eq:eta:d} (in a suitable way), then it must satisfy 
$$-1+1/n\le u_n(t,x) \le 1-1/n\quad \text{a.e.~in~}(0,T)\times \Omega,$$ otherwise
the constraint \eqref{eq:eta:d} could not be satisfied. For this reason,
for any fixed $n$, one can solve system \eqref{eq:mu:d}--\eqref{eq:eta:d}
by replacing the singular functions $\beta$ and $g$ outside the interval $[-1+1/2n,1-1/2n]$ with smooth extensions defined on the whole real line $\RR$. As a consequence,
the singularities due to $\beta$ and $g$ at $\pm 1$ simply disappear when dealing
with the system \eqref{eq:mu:d}--\eqref{eq:eta:d}. On the other hand, the price to pay for this is, of
course, the presence of the additional term $\zeta_n\in \Gamma_n(u_n)$. However, we shall see that the difficulty induced by $\Gamma_n$ 
is actually simpler to deal with comparing with those involving $\beta$, because it acts only 
on $u_n$, whereas the terms depending on $\beta$ may involve derivatives up to the second order in space.

Taking the above considerations into account, a well-posedness result for the approximate system \eqref{eq:mu:d}--\eqref{eq:eta:d} can be obtained 
by applying a similar argument for \cite[Theorem~3.1]{SP0}, where a six-order Cahn--Hilliard equation with nonlinear diffusions was studied. More precisely, we have
\bete\label{thm:dd}
 Assume that the hypotheses of\/ {\rm Theorem~\ref{thm:main}} are satisfied. 
 For any integer $n\geq 3$, let us assume in addition
 \begin{equation}\label{init:dd}
   -1 + 2n^{-1} \le u_0 \le 1 - 2n^{-1}
    \quext{a.e.~in }\,\Omega.
 \end{equation}
 Then there exists a unique function $u_n$, some function $\mu_n$,
 and some functional $\zeta_n$, with
 \begin{align}\label{rego:udd}
   & u_n\in H^1(0,T;V') \cap L^\infty(0,T;W)
    \cap L^4(0,T;H^3(\Omega)), \\
  \label{rego:udd2}
   & -1+n^{-1} \le u_n \le 1-n^{-1}
    \quext{a.e.~in }\,(0,T)\times \Omega,\\
  \label{rego:mudd}
   & \mu_n\in L^2(0,T;V), \\
  \label{rego:etadd}
   & \zeta_n\in L^2(0,T;V'),
 \end{align}
 satisfying\/ \eqref{eq:mu:d}--\eqref{eq:eta:d} for a.a.~$t\in(0,T)$, together with the initial condition
 \begin{equation}\label{E:initn}
   u_n|_{t=0} = u_0 \quext{a.e.~in }\,\Omega.
 \end{equation}
\ente
We remark that, in \cite[Theorem~3.1]{SP0}, it is only stated that
$u\in L^2(0,T;H^3(\Omega))$ in place of $L^4(0,T;H^3(\Omega))$
(cf.~the last of \eqref{rego:udd}). However, refining a bit
the estimates (like in Section \ref{sec:apriori}), we are able to show that the time-regularity
exponent for $u_n$ can be actually improved up to $4$.

It is also worth mentioning that the results obtained in Theorem \ref{thm:dd} can be improved
from the point of view of regularity, provided that the initial datum is smoother.
This case was just mentioned in \cite[Section~6.2]{SP0}, but not explicitly treated. 
However, proceeding as in Section~\ref{sec:long} below (for what
concerns the parabolic smoothing estimates), one can easily realize
that, supposing additionally
\begin{equation}\label{init:dd2}
   u_0 \in D(A^2) \quad \text{and} \quad A^2 u_0 \in V,
\end{equation}
the solution $(u_n,\mu_n,\zeta_n)$ to problem \eqref{eq:mu:d}--\eqref{eq:eta:d} with \eqref{E:initn}  
satisfies the additional regularity properties
\begin{align}\label{rego:udd+}
  & u_n\in W^{1,\infty}(0,T;V') \cap H^1(0,T;W)
    \cap L^\infty(0,T;H^3(\Omega)), \\
 \label{rego:mudd+}
  & \mu_n\in L^\infty(0,T;V),\\
 \label{rego:etadd+}
  & \zeta_n\in L^\infty(0,T;V').
\end{align}
We shall omit the proof of this fact here but just discuss ``informally'' the necessity of condition \eqref{init:dd2}. Actually, the point is that, to have
\eqref{rego:mudd+} starting from the initial time, it is needed
to know that $\mu_n|_{t=0}$ lies in $V$. Of course $\mu_n$ is just
an auxiliary variable and its regularity at $t=0$
should be deduced from
that of $u_0$ by comparison of terms in \eqref{eq:u:d} ``evaluated''
at the initial time. In view of the presence of the bi-Laplacian, this
gives back the condition \eqref{init:dd2}. In the same spirit, one may
observe that ``evaluating'' \eqref{eq:u:d} at $t=0$ one should
also face the term $\zeta_n|_{t=0}$. This is however just $0$ in
view of the fact that the support of the initial datum $u_0$ 
has been supposed to be strictly smaller than the domain
of $\mathcal{I}_n$ (cf.~\eqref{defId} and \eqref{init:dd}) in order to avoid possible concentration phenomena.

\smallskip

Now, in order to apply Theorem~\ref{thm:dd} and make use of the subsequent
observations mentioned above, we also need to regularize the initial datum. Namely, given $u_0$ as in \eqref{hp:init}, we have to construct a sequence of 
$u_{0,n}$ complying with the constraints \eqref{init:dd} and \eqref{init:dd2}
and additionally satisfying
\begin{equation}\label{init:dd4}
  u_{0,n} \to u_0 \quext{strongly in }\,W, \qquad
   \| \beta(u\zzn) \| \le C\big( 1 + \| \beta(u_0) \|\big),
\end{equation}
for some constant $C$ that is independent of $n$.
The construction of $u_{0,n}$ turns out to be a bit technical. First of all, we take
\begin{equation}\label{choice:u0n}
  v_{0,n}^{(1)}(x):= \Big(1-\frac2n\Big)u_0(x),\quad \forall\, n\in \mathbb{N},\ n\geq 3.
\end{equation}
Next, we define $v_{0,n}^{(2)}$ as the solution to the elliptic
problem (recalling that $A$ is the minus Neumann Laplacian)
\begin{equation}\label{ell:11}
  v_{0,n}^{(2)} + n^{-1} A v_{0,n}^{(2)} = v_{0,n}^{(1)}
\end{equation}
and then we iterate the procedure by putting
\begin{equation}\label{ell:12}
  u_{0,n} + n^{-1} A u_{0,n} = v_{0,n}^{(2)}.
\end{equation}
In this way, $v_{0,n}^{(1)}$ takes values in the interval $[-1+2/n,1-2/n]$
by construction and the same holds for both $v_{0,n}^{(2)}$ and $u_{0,n}$
thanks to the maximum principle. Hence, the constraint \eqref{init:dd} is fulfilled. 
Moreover, by the elliptic regularity theory one can check that the functions $u_{0,n}$ also satisfy the regularity requirement as in \eqref{init:dd2}. Finally, we prove \eqref{init:dd4}. Indeed, the first property can be checked in a straightforward way. To prove the second one,
we first observe that $\| \beta(v\zzn^{(1)}) \| \le \| \beta(u_0) \|$ by construction. Then, we consider a smooth convex function $\psi$ over $(-1,1)$ that explodes at $\pm1$ as fast as $\beta^2$. Testing (formally) \eqref{ell:11} by $\psi'(v_{0,n}^{(2)})$
and applying the convexity of $\psi$ (i.e., $\psi''\geq 0$), we can deduce that
\begin{equation}\label{ell:11b}
  \io \psi\big(v_{0,n}^{(2)}\big) \,\dix \le \io \psi\big(v_{0,n}^{(1)}\big) \,\dix.
\end{equation}
Repeating the same procedure for \eqref{ell:12} then one can easily arrive at the second property of \eqref{init:dd4}.

%As a consequence (see, e.g., \cite[Thm.~3.66]{At},
%the operators $\Gamma_n$ converge to $\Gamma$
%(the subdifferential of $I$ in the duality between $V$ and $V'$)
%in the sense of graph convergence in $V\times V'$.
%(THIS IS PROBABLY NOT USED, BECAUSE ALL TERMS ARE MANAGED
%DIRECTLY, SEE BELOW)
%
%\smallskip
%
%%%%%%%%%%%%%%%%%%%%%%%%%%%%%%%%%%%%%%%%%%%%%%%%%%%%%%%%%%%%%%%%%%555
\subsection{Justification of \textit{a priori} estimates via approximate solutions}
\label{subsec:justi}

Below we show that all the formal \textit{a priori} estimates performed
in Section \ref{sec:apriori} become rigorous once one considers
the solutions $(u_n,\zeta_n)$ to the approximate system \eqref{eq:mu:d}--\eqref{eq:eta:d} with regularized initial data $u_{0,n}$ constructed above. 
Indeed, the presence of the singular constraint $\Gamma_n$ automatically guarantees the separation property \eqref{rego:udd2}. 
As a consequence, we can treat $\beta$ and all its derivatives
as if they were smooth and bounded functions.
On the other hand, we should notice that the presence of $\Gamma_n$ gives rise to the occurrence
of some additional terms in the \textit{a priori} estimates. In order to show
that in fact all these new terms can be managed, we now revisit 
each \textit{a priori} bounds in the frame of this approximation.

\smallskip

\noindent%
{\bf Energy estimate.}~~%
When testing \eqref{eq:u:d} by $\partial_t u_{n}$ one has to deal with the term $\zeta_n$. We claim that
\begin{equation}\label{estn:11}
  \duav{\zeta_n, \partial_t u_{n}}_{V',V} 
   = \ddt \mathcal{I}_n(u_n),
\end{equation}
whence integration in time yields
\begin{equation}\label{estn:12}
  \itt \duav{\zeta_n, \partial_t u_{n}}_{V',V} \,\dis
   = \mathcal{I}_n(u_n(t)) - \mathcal{I}_n(u\zzn). 
\end{equation}
The first term on the \rhs\ of \eqref{estn:12} is nonnegative (so that
in particular the additional constraint $-1+1/n \le u_n \le 1-1/n$
keeps holding in time), while the second one is $0$ in
view of the assumption \eqref{init:dd}. Hence, once \eqref{estn:11} is
established, no further problems arise. To prove \eqref{estn:11},
we first notice that its \lhs\ makes sense because
$\partial_t u_{n}\in \LDV$ thanks to \eqref{rego:udd+}.
Then, noting as $\calR : V \to V'$ the Riesz operator,
we have
\begin{equation}\label{estn:13}
  \duav{\zeta_n, \partial_t u_{n}}_{V',V} 
   = ( \calR^{-1} \zeta_n, \partial_t u_{n} )_V
   = \ddt \mathcal{I}_n(u_n),
\end{equation}
where $( \cdot, \cdot )_V$ denotes the scalar product of $V$
and the second equality follows from the classical chain rule formula
for maximal monotone operators (see \cite[Lemma~3.3, p.~73]{Br}), 
once one observes that $\calR^{-1}\circ \Gamma_n$ indeed coincides with the subdifferential of the functional $\mathcal{I}_n$ with respect to the scalar product of $V$.

\smallskip

\noindent%
{\bf Second estimate.}~~%
Here dealing with the additional term $\zeta_n$ requires some
additional care. First of all, by \eqref{choice:u0n}--\eqref{ell:12},
one can easily check that
\begin{equation}\label{estn:20a}
  \overline{u\zzn} = \Big( 1 - \frac2n \Big) \overline{u_0}.
\end{equation}
Namely, the mean value of $u\zzn$ is closer to $0$ compared
to the mean value of $u_0$. Besides, this property keeps holding also for $t>0$ thanks to the mass conservation. As a consequence, the argument leading to \eqref{co:21} holds uniformly in $n$,
in particular, the constants $\kappa$ and $C$ in \eqref{co:21}
can be taken independent of $n$.

Next, we notice that a new term appears in the \lhs\ of the analogue
of~\eqref{co:22}, namely the duality $\duav{\zeta_n,u_n - \overline{u_n}}_{V',V}$.
To manage it, we need to recall that
(cf., e.g., \cite[Theorem~2.2]{SP0} and the related discussion),
for almost all $t\in (0,T)$, $\zeta_n$
can be interpreted as a measure $\nu_n$ that may be decomposed
as the sum of an absolutely continuous part
$\nu_{n,a}\in L^1(\Omega)$ and a part $\nu_{n,s}$
being singular with respect to the Lebesgue measure.
Moreover, $\nu_{n,a}(x)\in \de I_{[-1+1/n,1-1/n]}(u_n(x))$
for a.e.~$x\in\Omega$. In particular, $\nu_{n,a}$ is ``supported''
on the set where $|u_n|=1-1/n$ and has the same sign as $u_n$.
A similar property holds for the singular part (see \cite{SP0}
for details).  Hence, using the fact that $u_n$ is continuous up to the boundary of $\Omega$
thanks to \eqref{rego:udd} and the continuous embedding $W\subset C(\barO)$,
the new term can be controlled observing that
\begin{align}\no
  \duav{\zeta_n,u_n- \overline{u_n}}_{V',V}
    & = \io \nu_{n,a} (u_n-\overline{u_n}) \dix
     + \ibaro (u_n-\overline{u_n}) \deriv\! \nu_{n,s}\\
 \label{key:dd0}
    & \ge \kappa_0 \big( \| \nu_{n,a} \|_{L^1(\Omega)}
     + | \nu_{n,s} |(\barO) \big).
\end{align}
Here, $| \nu_{n,s} |$ denotes the total variation of
the measure $\nu_{n,s}$ and we have used the fact that,
in view of assumption \eqref{hp:mean} and of property 
\eqref{estn:20a}, there exists a constant $\kappa_0>0$ depending
only on $m$ such that, at least for $n\in\NN$ sufficiently large,
it holds
\begin{equation}\label{co:41:da}
  - (1 - 1/n) + \kappa_0 \leq  \overline{u_n} \leq (1 - 1/n) - \kappa_0.
\end{equation}
For instance, here we can take $\kappa_0=\frac12 \min\{1-m,\, 1+m\}$.
As a result, the new contribution in \eqref{key:dd0}
is positive. Indeed, we see from 
the analogue of \eqref{co:22} for $u_n$ that the above new term provides the additional information
\begin{equation}\label{co:41:dd0}
  \big( \| \nu_{n,a} \|_{L^1(\Omega)}
     + | \nu_{n,s} |(\barO) \big)
     \le \Lambda_1(t),
\end{equation}
with certain function $\Lambda_1\in L^2(0,T)$. As a further consequence, when we integrate
\eqref{eq:uom1} at the level $n$, on the \rhs\ of \eqref{st:23a}
appears the new term $|\Omega|^{-1}\duav{\zeta_n,1}_{V',V}$, which in principle
has no sign and needs to be controlled. However, it is clear that
\begin{equation}\label{co:41:dd1}
  |\duav{\zeta_n,1}_{V',V}| \le \big( \| \nu_{n,a} \|_{L^1(\Omega)}
     + | \nu_{n,s} |(\barO) \big)
\end{equation}
for almost all $t\in (0,T)$. Therefore, thanks to \eqref{co:41:dd0}, \eqref{st:24a}
and \eqref{st:24} keep holding (in the sense that they provide uniform \textit{a priori}
estimates  with respect to~$n$).

\smallskip

\noindent%
{\bf Third estimate.}~~%
We have for almost all $t\in (0,T)$,
\begin{equation}\label{estn:31}
  \duav{\zeta_n, Au_n}_{V',V} \ge 0,
\end{equation}
thanks to the regularity properties \eqref{rego:udd}, \eqref{rego:etadd},
 as well as the result \cite[Lemma~2.4]{SP0}. Thus, \eqref{st:31} holds uniformly with respect to $n$.

\smallskip

\noindent%
{\bf Fourth estimate.}~~%
First of all, we notice that the
function $\beta$ is smooth and bounded with
all its derivatives in the interval $[-1+1/n,1-1/n]$
(which includes the range of $u_n$), hence
$\beta(u_n)$ has the same regularity of $u_n$ and it can 
be used as a test function. So, the estimate therein can still be
performed with the additional term
\begin{equation}\label{co:40:dd}
  \duav{\zeta_n,\beta(u_n)}_{V',V}
\end{equation}
to be handled. Now, since $\beta(u_n)$ is continuous
up to the boundary of $\Omega$ (this follows from \eqref{rego:udd}
and the continuous embedding $W\subset C(\barO)$), proceeding
as above, we infer that 
\begin{align}\no
  \duav{\zeta_n,\beta(u_n)}_{V',V}
    & = \io \nu_{n,a} \beta(u_n) \dix
     + \ibaro \beta(u_n) \deriv\! \nu_{n,s}\\
 \label{key:dd}
    & = \beta(1-1/n) \big(
     \| \nu_{n,a} \|_{L^1(\Omega)}
     + | \nu_{n,s} |(\barO) \big),
\end{align}
where $| \nu_{n,s} |$ denotes the total variation of
the measure $\nu_{n,s}$. We note that this contribution
is positive, and as a consequence, \eqref{co:41b} holds uniformly with respect to $n$ (thus \eqref{st:41}).

\smallskip

\noindent%
{\bf Handling the extra singular term $\zeta_n$ as $n\nearrow +\infty$}.~~%
Finally, we show that the estimate \eqref{key:dd} also helps us to get rid of the additional term
$\zeta_n$ in the limit $n\nearrow +\infty$. Indeed, with the notation
of~\eqref{co:41b}, we now have the additional information
\begin{equation}\label{co:41:dd}
  \beta(1-1/n) \big( \| \nu_{n,a} \|_{L^1(\Omega)}
     + | \nu_{n,s} |(\barO) \big)
     \le \Lambda(t),
\end{equation}
with $\Lambda\in L^2(0,T)$ (cf.~\eqref{co:41c}).
Hence, squaring \eqref{co:41:dd} and integrating in time,
we obtain
\begin{equation}\label{co:42:dd}
   \| \nu_{n,a} \|_{L^2(0,T;L^1(\Omega))}^2
     + \iTT | \nu_{n,s} |^2(\barO) \, \dit
     \le \frac{C}{\beta^2(1-1/n)},
\end{equation}
where the constant $C$ is independent of $n$. Then, as one lets $n\nearrow +\infty$, from the fact
$\beta(1-1/n) \nearrow +\infty$ it follows that 
\begin{equation}\label{co:43:dd}
   \nu_{n} \to 0 \quext{strongly in }\, L^2(0,T;\calM(\barO)),
\end{equation}
where $\calM(\barO)=C(\barO)'$ is the space of (signed) Borel measures
on $\barO$. In particular, for all $\fhi\in V\cap C(\barO)$, we have
\begin{equation}\label{co:44:dd}
   \duav{\zeta_n,\fhi}_{V',V} =\io \fhi\,\deriv\!\nu_{n} \to 0
     \quext{strongly in }\, L^2(0,T),
\end{equation}
whence, in view of the continuous embedding $W\subset V\cap C(\barO)$, we obtain
(at least)
\begin{equation}\label{co:45:dd}
   \zeta_n \to 0 \quext{weakly in }\, L^2(0,T;W').
\end{equation}
Namely, the extra singular term $\zeta_n$ in the approximate system \eqref{eq:mu:d}--\eqref{eq:eta:d} disappears in the limit $n\nearrow +\infty$, as desired.

%%%%%%%%%%%%%%%%%%%%%%%%%%%%%%%%%%%%%%%%%%%%%%%%%%%%%%%%%%%%%%%%%%%%%%%%%%%%%%%%%%%%%%%%%%%%%

\section{Uniqueness}
\label{sec:uniq}

In this section, we will provide an alternative formulation of problem
\eqref{eq:mu}--\eqref{ini} ruled by an abstract operator of subdifferential type.
This formulation will be {\sl weaker}\/ with respect to that
provided in Definition \ref{def:weak}. In particular, we
shall prove that the solutions given by Theorem~\ref{thm:main}
also solve this subdifferential formulation. Since uniqueness
for the subdifferential formulation can be proved by standard monotone operator tools,
then uniqueness will also hold for the weak solutions
in the sense of Definition \ref{def:weak}.

\subsection{Subdifferential interpretation}
\label{sec:subinter}

We start with introducing a number of preliminaries. The following simple
property is proved just for the reader's convenience:
\bele\label{lem:trunc}
 Let $v\in V$ such that $-1 \le v(x) \le 1$ a.e.~in $\Omega$.
 We define the truncation
 \begin{equation}\label{trunc}
   v_n = T_n v = \max\big\{ - 1 + 1/n , \min\{ v , 1 - 1/n \} \big\}.
 \end{equation}
 Then, $v_n \to v$ uniformly and strongly in $V$ as $n\nearrow +\infty$.
\enle
\begin{proof}
Convergence in $L^2$ and uniform convergence are obvious. To prove
convergence in~$V$, we notice that
\begin{equation}\label{trunc:11}
  \| \nabla v_n - \nabla v \|^2
   = \int_{\{ |v|\ge 1-1/n \} } | \nabla v|^2 \,\dix
   \to \int_{\{ |v| = 1 \} } | \nabla v|^2 \,\dix
   = 0.
\end{equation}
In the above expression we have used the chain rule formula for
Sobolev functions $\nabla (G\circ v) = G'(v) \nabla v$ holding for $v\in V$ and
Lipschitz operator $G$ (here applied with $G=T_n$) together with
Lebesgue's dominated convergence theorem.
\end{proof}
\noindent%
The chain rule for Sobolev functions plays an important role in the
above proof. In the sequel we shall need a more refined version of it,
which is stated below for the reader's convenience
in a form suitable for our purposes. For the proof one can refer, e.g.,
to \cite[Theorem~2.1]{MaMi}, where a more general statement is given.
\bele\label{lem:chain}
 Let $\gamma:\RR \to \RR$ be an absolutely continuous function
 (hence, in particular, let $\gamma'\in L^1(\RR)$). Let $v\in V$
 and let us assume that $\gamma'(v)=\gamma'\circ v \in L^2(\Omega)$.
 Then $\gamma(v)=\gamma\circ v\in W^{1,1}(\Omega)$ and
 $\partial_{x_i} \gamma(v) = \gamma'(v) \partial_{x_i} v$
 for $i=1,...,3$, where the product on the
 \rhs\ is intended to be $0$ whenever $\partial_{x_i} v=0$.
\enle
Let now $u$ be a weak solution in the sense of Definition~\ref{def:weak}
and let us set for later convenience 
$$
\gamma(r):=\arcsin r\quad \text{and}\quad b(r)=\gamma'(r)=\frac{1}{(1-r^2)^{1/2}}.
$$ 
Here $\gamma$ is defined only for $r\in [-1,1]$
but it is clear that it can be extended to the whole real line $\RR$ 
in such a way that its extension lies in $W^{1,1}(\RR)$.
Then, from \eqref{reg:bu2} we have at least
$b(u) \in L^4(0,T;L^2(\Omega))$. Moreover, in view of \eqref{reg:u}, we have 
$\nabla u\in L^4(0,T;H^2(\Omega))\subset L^4(0,T;L^\infty(\Omega))$.
Hence, we deduce that 
\begin{align*}
 & \nabla \gamma(u) = b(u)\nabla u
   = \frac1{(1-u^2)^{1/2}} \nabla u~~\text{a.e.~in }\,(0,T),
\end{align*}
and 
\begin{align*}
& \gamma(u)=\arcsin u\in L^2(0,T;V).
\end{align*}
We can now define the functional $\mathcal{J}:V \to [0,+\infty]$ which is at the
core of the subdifferential formulation of our problem:
\begin{equation}\label{defi:J}
  \mathcal{J}(v):=
   \begin{cases}\displaystyle
     \io | \nabla (\arcsin v) |^2\,\dix
    & \text{~~if }\, \arcsin v \in V,\\
    + \infty & \text{~~otherwise}.
   \end{cases}
\end{equation}
Here, we are implicitly asking that the domain of
the functional $\mathcal{J}$, i.e., the set $D(\mathcal{J})$ where it takes finite values,
may only contain those functions $v\in V$ such that
$-1\le v\le 1$ a.e.~in~$\Omega$ and $\arcsin v \in V$.
On the other hand, for $v\in D(\mathcal{J})$
the set where $|v|=1$ may be large. In particular,
the constant functions $v\equiv\pm 1$ lie in $D(\mathcal{J})$.
Besides, we note also that, by the above argument, if $u$ is a
weak solution to problem \eqref{eq:mu}--\eqref{ini} as in Definition~\ref{def:weak}, then $u(t)\in D(\mathcal{J})$ for a.e.~$t\in (0,T)$. 

A number of additional properties of $\mathcal{J}$ are summarized in the following lemma.
\bele\label{lemma:J}
 The functional $\mathcal{J}$ is convex and lower semicontinuous on $V$.
 Moreover, if $v\in V$ satisfies $-1 + \epsi \le v(x) \le 1-\epsi$
 for almost every $x\in \Omega$ and some $\epsi\in(0,1)$,
 then $\mathcal{J}$ is G\^ateaux-differentiable at the point $v$, with
 respect to the norm of $V\cap L^\infty(\Omega)$ and
 \begin{equation}\label{J:gat}
    \duav{ D\mathcal{J}(v),\fhi }
     = \io a(v) \nabla v \cdot \nabla \fhi \,\dix
      + \frac12 \io a'(v) | \nabla v |^2 \fhi\,\dix
      \qquad \perogni \fhi \in V \cap L^\infty(\Omega),
 \end{equation}
where the function $a(\cdot)$ is defined as in \eqref{defia}.
\enle
\begin{proof}
Let us start with showing that $\mathcal{J}$ is lower semicontinuous.
To this aim, let $\{v_n\}\subset V$ with $v_n \to v$ in $V$.
Then, we can suppose without loss of generality that the sequence
$\{\mathcal{J}(v_n)\}$ is bounded and that $v_n$ as well as $\nabla v_n$ tend respectively
to $v$ and $\nabla v$ pointwisely (in fact this holds
at least for a subsequence). This implies in particular
that both $v_n$ and $v$ take values in the interval $[-1,1]$.
As a consequence, it holds $\arcsin v_n \to \arcsin v$ pointwisely.
Since the function $\arcsin$ is bounded on $[-1,1]$,
we deduce that $\arcsin v_n \to \arcsin v$ weakly
in $H$. Consequently, we have $\nabla (\arcsin v_n) \to \nabla (\arcsin v)$
in the sense of distributions. Moreover, from
the boundedness of $\{\mathcal{J}(v_n)\}$ we further deduce that
$\nabla (\arcsin v_n) \to \nabla (\arcsin v)$ weakly in $H$.
From the semicontinuity of the $H$-norm with respect to
weak convergence, we then infer
\begin{equation}\label{semico}
  \mathcal{J}(v) = \| \nabla (\arcsin v) \|^2
   \le \liminf_{n\nearrow +\infty} \| \nabla (\arcsin v_n) \|^2
   = \liminf_{n\nearrow +\infty} \mathcal{J}(v_n),
\end{equation}
as desired.

Next, for $v\in V$ satisfying $-1 + \epsi \le v \le 1-\epsi$
a.e.~in~$\Omega$, let us show the G\^ateaux-differentiability of $\mathcal{J}$ at~$v$. 
Let $\fhi\in V\cap L^\infty(\Omega)$ and $t\in \RR$.
Then, for $|t|$ small enough, we have $|v(x) + t \fhi(x)| \le 1-\epsi/2$
a.e.~in~$\Omega$. In view of the fact that the function $\arcsin$
is smooth and bounded with all its derivatives
in the interval $[-1+\epsi/2,1-\epsi/2]$, we may apply to the functions
$\arcsin v$, $\arcsin \fhi$ and $\arcsin(v+t\fhi)$ the (standard)
chain rule formula of Sobolev spaces. Namely, we have
$$
  J(v) = \io | \nabla (\arcsin v) |^2\,\dix
   = \io \frac{a(v)}2 | \nabla v |^2\,\dix
$$
with similar equalities holding for $\fhi$ and $v+t\fhi$. Here, the function $a(\cdot)$ is defined as in \eqref{defia}. 
Moreover, one can easily check that, as $t\to 0$,
\begin{align}\no
   & \frac{\mathcal{J}(v+t\fhi)-\mathcal{J}(v)}t\\
   \no
    & \quad = \frac12 \io \frac{a(v+t\fhi)-a(v)}t | \nabla (v+t\fhi) |^2 \,\dix
      + \frac12 \io a(v) \frac{| \nabla (v+t\fhi) |^2 - |\nabla v|^2}t \,\dix\\
 \label{J:gat3}
  &\quad \to \frac12 \io a'(v)  | \nabla v |^2 \fhi\,\dix
     + \io a(v) \nabla v \cdot \nabla \fhi \,\dix,
\end{align}
where we notice that, as $v$ is an assigned function as above, it holds
\begin{equation}\label{J:gat4}
  \bigg | \frac12 \io a'(v) | \nabla v |^2 \fhi \,\dix
     + \io a(v) \nabla v \cdot \nabla \fhi \,\dix \bigg|
    \le C(v) \big( \| \fhi \|_V + \| \fhi \|_{L^\infty(\Omega)} \big).
\end{equation}
Namely, the G\^ateaux derivative $D\mathcal{J}(v)$ exists and it acts as a bounded linear
functional on $V\cap L^\infty(\Omega)$, as desired.

In order to show the convexity of $\mathcal{J}$, let us
first take a couple of functions $u,v\in V$
both taking values in the interval $[-1+\epsi,1-\epsi]$
for some $\epsi\in (0,1)$. Let us also set
$\Psi(t):=\mathcal{J}(v+t(u-v))$ for $t\in[0,1]$. Then, using again the fact that
the function $\arcsin$ is bounded with all its derivatives
once its argument remains inside the interval $[-1+\epsi,1-\epsi]$,
we can compute directly the second derivative $\Psi''(t)$.
In particular, following the lines of \cite[Theorem~6.1]{SP0},
one can show by some lengthy but otherwise elementary computations
that $\Psi''(t)\ge 0$ for any $t\in[0,1]$. This clearly
implies that $\Psi$ is convex. As a consequence, we get 
\begin{equation}\label{J:convx}
  \mathcal{J}(t u + (1-t) v) = \Psi(t)
   \le t \Psi(1) + (1-t) \Psi(0)
    = t \mathcal{J}(u) + (1-t) \mathcal{J}(v),
\end{equation}
for any $t\in[0,1]$, i.e., $\mathcal{J}$ is convex.

Next, let us take any couple $u,v\in D(\mathcal{J})$. 
We set $u_n:=T_n u$ and $v_n:=T_n v$ (cf.~\eqref{trunc}).
Let us also observe that, as $n\nearrow +\infty$,
\begin{align}\no
  \mathcal{J}(u_n) 
    & = \io | \nabla (\arcsin u_n) |^2\,\dix
    = \int_{\{|u|\le 1-n^{-1} \} } | \nabla (\arcsin u) |^2\,\dix \\
    \no
    & = \int_{ \{ |\arcsin u|\le \arcsin(1-n^{-1}) \} } | \nabla (\arcsin u) |^2\,\dix\\
    \label{J:gat4b}
    &
    \to \io | \nabla (\arcsin u) |^2\,\dix = \mathcal{J}(u),
\end{align}
where we have applied Lemma~\ref{lem:trunc} to the function
$\arcsin u$. Since both $u_n$ and $v_n$ take values in
the interval $[-1+1/n,1-1/n]$, we have
\begin{equation}\label{J:gat4x}
  \mathcal{J}(t u_n + (1-t) v_n)
   \le t \mathcal{J}(u_n) + (1-t) \mathcal{J}(v_n).
\end{equation}
for all $t\in[0,1]$. Then, in view of the fact that
$t u_n + (1-t) v_n$ tends to $t u + (1-t) v$ in $V$,
by lower semicontinuity we obtain that
\begin{align}\no
   \mathcal{J}(t u + (1-t) v)
   &  \le \liminf_{n\nearrow +\infty} \mathcal{J}(t u_n + (1-t) v_n)\no\\
   &  \le \liminf_{n\nearrow +\infty} \big( t \mathcal{J}(u_n) + (1-t) \mathcal{J}(v_n) \big) \no\\
   &  = t \lim_{n\nearrow +\infty} \mathcal{J}(u_n) + (1-t) \lim_{n\nearrow +\infty} \mathcal{J}(v_n)\no\\
   &  = t \mathcal{J}(u) + (1-t) \mathcal{J}(v),\quad \forall\, t\in [0,1].
   \label{J:gat4y}
\end{align}
Since the above inequality holds for any couple $u,v\in D(\mathcal{J})$, we have proved the convexity of $\mathcal{J}$. This concludes the proof of the lemma.
\end{proof}
\beos\label{twofunct}
 One may also introduce the (related) functional
 $\mathcal{K}:V\to [0,+\infty]$ defined by
 \begin{equation}\label{defi:K}
  \mathcal{K}(v):=
    \begin{cases}\displaystyle
      \io \frac{a(v)}2 | \nabla v |^2\,\dix
     & \text{~~if }\,a(v) |\nabla v|^2 \in L^1(\Omega),\\
     + \infty & \text{~~otherwise},
    \end{cases}
 \end{equation}
 where it is intended that the domain $D(\mathcal{K})$ of $\mathcal{K}$ is given by those functions
 $v\in V$ satisfying $-1 < v < 1$ a.e.~in $\Omega$ and such that the above integral is finite. However, in this case we cannot admit $|v|=1$ on a set of strictly positive measure, otherwise it would not be clear how to interpret the integrand due to the definition of $a(\cdot)$ (see \eqref{defia}).
 In particular, it is worth observing that the functionals $\mathcal{K}$ and $\mathcal{J}$ do not coincide. 
 More precisely, using Lemma~\ref{lem:chain}, one can prove that
 $D(\mathcal{K})\subset D(\mathcal{J})$ and that $\mathcal{K}$ and $\mathcal{J}$ only coincide with each other on $D(\mathcal{K})$. On the other
 hand, $D(\mathcal{J})$ is strictly larger than $D(\mathcal{K})$ (for instance the constant
 function $u\equiv1$ belongs to $D(\mathcal{J})$ but not to $D(\mathcal{K})$). For the same reason, $\mathcal{K}$ is not lower semicontinuous on $V$ (consider the sequence  $\{u_n\}$ with $u_n\equiv 1-n^{-1}$).
\eddos
Using the functional $\mathcal{J}$, we can finally introduce the desired abstract reformulation of the original problem
\eqref{eq:mu}--\eqref{ini} in terms of a subdifferential operator.
To this aim, we first observe that the restriction of $\mathcal{J}$ to the
space $W$ (still indicated by $\mathcal{J}$ for simplicity) is convex
and lower semicontinuous on $W$. Hence, we can denote by
$\calB$ the subdifferential of $\mathcal{J}$ with respect
to the duality pairing between $W$ and $W'$. Namely,
for $\xi\in W'$, $v\in W$, we set
\begin{equation}\label{de:J}
  \xi \in \calB(v)~~\text{if and only if}~~
   \duav{\xi,z-v}_{W',W} + \mathcal{J}(v) \le \mathcal{J}(z)\quad\perogni z\in W.
\end{equation}
Thanks to standard results on subdifferential operators in Hilbert
spaces (see e.g.,\cite{Ba,Br}), $\calB$ is a maximal monotone, possibly
multivalued, operator from $W$ to $2^{W'}$. Using the operator $\calB$,
we can define \textit{sudifferential solutions} as follows:
\bede\label{def:subd}
 A triple $(u,\mu,\zeta)$ is called a subdifferential solution to problem
\eqref{eq:mu}--\eqref{ini} over the time interval $(0,T)$, provided that:\\[2mm]
 (A)~~The regularity conditions~\eqref{reg:u}, \eqref{reg:bu2}, \eqref{reg:mu}
 are satisfied, together with
 \begin{equation} \label{reg:subd}
   \zeta \in L^2(0,T;W').
 \end{equation}
 (B)~~The following weak counterparts of equations\/ \eqref{E:mu} and \eqref{E:u}
 are satisfied for a.a. $t\in(0,T)$:
 \begin{align}\label{E:mus}
  & \partial_t u + A \mu = 0 \quext{in }\,V',\\
  \label{E:us}
   & \mu = A^2 u + \zeta
    + \beta(u) \beta'(u)
    - (2\lambda - \eta) A u + g(u)
     \quext{in }\,W',\\
  \label{E:zetas}
   & \zeta \in \calB(u) \quext{in }\,W'.
 \end{align}
 (C)~~The initial condition is satisfied in the sense of~\eqref{E:init}.
\edde

Next, we establish the relation between the weak solution given by Definition \ref{def:weak} and the  sudifferential solution given by Definition \ref{def:subd}. 
\bele\label{lemma:subd}
 Let $(u, \mu)$ be a weak solution in the sense of\/
 {\rm Definition~\ref{def:weak}} on $(0,T)$.  For almost all $t\in(0,T)$,
 let us also  define
 \begin{equation}\label{defi:zeta}
   \zeta := 2 A \beta(u) + \beta''(u) | \nabla u |^2
     = \mu - A^2 u - \beta(u)\beta'(u) + ( 2\lambda - \eta ) Au - g(u).
 \end{equation}
 Then,
 \begin{equation}\label{zeta:subd}
   \zeta(t) \in \calB(u(t)) \quext{for almost all }\,t\in (0,T)
 \end{equation}
 and the triple $(u,\mu,\zeta)$ is a subdifferential solution
 in the sense of\/ {\rm Definition~\ref{def:subd}}.
\enle
\begin{proof}
By the regularity properties satisfied by weak solutions it is easy to
check that, at least,
\begin{equation}\label{rego:zeta}
  \zeta \in L^2(0,T;L^1(\Omega)) + L^2(0,T;V')
   \subset L^2(0,T;W'),
\end{equation}
the last inclusion following from the standard Sobolev
embedding theorem. Moreover, by \eqref{defi:zeta},
for a.a. $t\in (0,T)$, there holds
\begin{equation}\label{eq:subd:x}
  \duav{\zeta,\fhi}_{W',W}
   = 2 \io \nabla \beta(u) \cdot \nabla \fhi\dix
    + \frac12 \io a'(u) | \nabla u |^2 \fhi \dix
    \quad \perogni \fhi \in W.
\end{equation}
Then, verifying \eqref{zeta:subd} amounts to prove that,
for any $v\in D(\mathcal{J}) \cap W$ and a.a.~$t\in(0,T)$, it holds
\begin{equation}\label{subd:11}
  \duav{\zeta,v-u}_{W',W}
    \le \io | \nabla (\arcsin v) |^2\dix
     - \io | \nabla (\arcsin u) |^2\dix
    = \mathcal{J}(v) - \mathcal{J}(u).
\end{equation}
This will be shown by a truncation argument. Fix $t\in (0,T)$ and
let $u_n:= T_n u$ (cf. \eqref{trunc}). Set also $E:=V\cap L^\infty(\Omega)$
(which is a Banach space with the natural norm) and
let $\zeta_n$ be defined by
\begin{align}\no
  \duav{\zeta_n,\fhi}_{E',E}
   & := 2 \io \nabla \beta(u_n) \cdot \nabla \fhi \dix
    + \frac12 \io a'(u_n) | \nabla u_n |^2 \fhi \dix\\
 \label{subd:12}
   &\ \,= \io a(u_n) \nabla u_n \cdot \nabla \fhi \dix
    + \frac12 \io a'(u_n) | \nabla u_n |^2 \fhi \dix,
\end{align}
for any $\fhi\in E$. Then by a direct check we can verify
that $\zeta_n \in E'$. We also point out that the two expressions of the
above \rhs\ are equivalent because $-1+1/n\le u_n(x) \le 1-1/n$
for every $x\in \Omega$. Let now $\fhi\in W$ (keeping in mind
that $W\subset E$ continuously). We assume that $y\in \calB(u_n)$. Then, for any $z\in W$, $h\in\RR$,
using the definition of subdifferential, we get
\begin{equation}\label{subd:a1}
  \mathcal{J}(u_n + h z) \ge \mathcal{J}(u_n) + h \duav{y,z}_{W',W}.
\end{equation}
Consequently, we infer that
\begin{equation}\label{subd:a2}
  \frac{ \mathcal{J}(u_n + h z) - \mathcal{J}(u_n) }h \ge \duav{y,z}_{W',W}
   \quext{for all}\,h>0,
\end{equation}
with the opposite inequality holding for $h<0$.
Now, recalling Lemma~\ref{lemma:J}, $\mathcal{J}$ is G\^ateaux differentiable
at $u_n$ with respect to the norm of $E$ and its G\^ateaux derivative
coincides with $\zeta_n$. Hence, since $z\in W\subset E$,
taking the limit $h\searrow 0$ in \eqref{subd:a2}
and the limit $h\nearrow 0$ in its analogue for $h<0$, we
easily infer
\begin{equation}\label{subd:a3}
  \duav{\zeta_n,z}_{E',E} = \duav{y,z}_{W',W}
   \quext{for all }\,z\in W.
\end{equation}

Hence, in view of the density of $W$ in $E$
we have obtained that, if $y$ is an element of $\calB(u_n)$
(i.e., the $W$-subdifferential of $\mathcal{J}$ at $u_n$), then $y$ admits a unique
extension as a linear and continuous functional defined on $E$
and this extension coincides with $\zeta_n$, that is, the
G\^ateaux derivative of $\mathcal{J}$ at $u_n$. In other words,
with a small abuse of notation, we can write $\calB(u_n)=\{\zeta_n\}$.
In particular, $\calB(u_n)$ contains a single element. 
As a consequence, for any $v\in W$ we can write
\begin{align}\no
  & 2 \io \nabla \beta(u_n) \cdot \nabla ( v - u_n ) \, \dix
    + \frac12 \io a'(u_n) | \nabla u_n |^2  ( v - u_n ) \, \dix\\
\no
 &  \mbox{}~~~~~  =  \duav{\zeta_n,v-u_n}_{W',W} \\
 \label{subd:13}
  & \mbox{}~~~~~
   \le \mathcal{J}(v) - \mathcal{J}(u_n)\\
   \no
   &\mbox{}~~~~~
      = \io | \nabla (\arcsin v) |^2 \, \dix
     - \io | \nabla (\arcsin u_n) |^2 \, \dix,
\end{align}
where the \rhs\ is intended to be $+\infty$ in the case
when $v\not\in D(\mathcal{J})$.
Now, recalling that $u(\cdot)\in W$ for a.a. $t\in(0,T)$, we can
plug $\fhi=u$ in \eqref{eq:subd:x}. This implies that, for a.a. $t\in(0,T)$,
there holds
\begin{equation}\label{eq:subd2}
   2 \io \nabla \beta(u) \cdot \nabla u\,\dix
    + \frac12 \io a'(u) | \nabla u |^2 u \, \dix
   \le \| \zeta \|_{W'} \| u \|_{L^\infty(0,T;W)}
   \le C \| \zeta \|_{W'}.
\end{equation}
Hence, by \eqref{rego:zeta}, the integrals on the \lhs\
of \eqref{eq:subd2} are finite for a.a. $t\in(0,T)$.
Moreover, due to Lemma~\ref{lem:trunc},
we observe that $u_n(\cdot) \to u(\cdot)$ strongly in $V$
for a.a. $t\in (0,T)$. Besides, we have at least
$\beta(u_n(\cdot)) \to \beta(u(\cdot))$ weakly in $V$.
Indeed, $\beta(u(\cdot))\in V$ for a.a. $t\in(0,T)$ thanks
to \eqref{reg:bu1}. As a consequence, we obtain,
for a.a. $t\in(0,T)$ that
\begin{equation} \label{leb:13}
  2 \io \nabla \beta(u_n) \cdot \nabla ( v - u_n ) \, \dix
   \to 2 \io \nabla \beta(u) \cdot \nabla ( v - u ) \, \dix.
\end{equation}
On the other hand, observing that
$|u_n| \le |u|$ and $|\nabla u_n| \le |\nabla u|$
a.e.~in~$(0,T)\times \Omega$, we can apply Lebesgue's dominated
convergence theorem to obtain
\begin{equation} \label{leb:12}
  \io a'(u_n) | \nabla u_n |^2  u_n \, \dix
    \to \io a'(u) | \nabla u |^2  u \, \dix.
\end{equation}
Next, in view of the fact that $v\in W \subset L^\infty(\Omega)$, it also
follows from Lebesgue's theorem that
\begin{equation} \label{leb:14}
  \io a'(u_n) | \nabla u_n |^2  v \, \dix
    \to \io a'(u) | \nabla u |^2  v \, \dix.
\end{equation}
Let us now take the limit $n\nearrow \infty$ in~\eqref{subd:13}.
Applying \eqref{leb:13}--\eqref{leb:14}, and noting that 
$$
\arcsin u_n(\cdot) \to \arcsin u(\cdot)\quad \text{strongly in }V,
$$
thanks to Lemma~\ref{lem:trunc},
we then deduce that, for a.a. $t\in(0,T)$, there holds
\begin{align} 
  & 2 \io \nabla \beta(u) \cdot \nabla ( v - u ) \, \dix
    + \frac12 \io a'(u) | \nabla u |^2  ( v - u ) \, \dix
  \le \mathcal{J}(v) - \mathcal{J}(u).
  \label{subd:14}
\end{align}
Now, the \lhs\ coincides with $\duav{\zeta, v-u}_{W',W}$ thanks to~\eqref{eq:subd:x}. 
Hence, we have obtained the conclusion \eqref{subd:11}, which completes the proof.
\end{proof}

\beos\label{rem:subd}
 (1) Lemma \ref{lemma:subd} implies that any weak solution is a subdifferential solution. Thus, Theorem~\ref{thm:main} indeed provides the existence of a subdifferential solution to problem
\eqref{eq:mu}--\eqref{ini}. 
 
(2) Since the concept of subdifferential solution is weaker than that of weak solution, it may happen, at least in principle, 
 that a subdifferential solution exists under certain weaker conditions on the initial datum $u_0\in W$. This is however not 
 expected in view of the fact that the finiteness of the initial energy $\mathcal{E}(u_0)$ also implies the second condition in \eqref{hp:init}. 
 
(3) A full characterization of the elements of the abstract operator $\calB(u)$ may be rather complicated. On the other hand, for weak solutions in the sense of Definition~\ref{def:weak}, the nonlinear diffusion part of equation \eqref{E:u} can be regarded as an element of $\calB(u)$ satisfying additional regularity properties that permit us to interpret it
 in a ``pointwise'' sense.
\eddos

%%%%%%%%%%%%%%%%%%%%%%%%%%%%%%%%%%%%%%%%%%%%%%%%%%%%%%%%%%%%%%%%%%%%%%%%%%%%%%%%%%%%%%%%%%%%555

\subsection{Proof of Theorem \ref{thm:main}: the uniqueness part}
\label{subsec:end}

In view of Lemma \ref{lemma:subd}, if we can prove the uniqueness of subdifferential solutions, then we immediately obtain that, from any initial datum $u_0$ satisfying \eqref{hp:init}, emanates one and only one weak solution to problem \eqref{eq:mu}--\eqref{ini}.

To this end, we derive a continuous dependence estimate for two subdifferential solutions $(u_i, \mu_i, \zeta_i)$, $i=1,2$. By definition, we have for a.a. $t\in(0,T)$,
\begin{equation}\label{uniq:11}
  \mu_i = A^2 u_i + \zeta_i + \beta(u_i)\beta'(u_i) - ( 2\lambda - \eta ) A u_i + g(u_i),
\end{equation}
with
\begin{equation}\label{uniq:12}
  \zeta_i \in \calB(u_i).
\end{equation}
Note that if $u_i$ are weak solutions, we further have 
\begin{equation}\no
  \zeta_i = 2 A \beta(u_i) + \beta''(u_i) | \nabla u_i |^2
   \in \calB(u_i) \quext{a.e.~in }\, (0,T).
\end{equation}
Let us set $u:=u_1 - u_2$ and $\mu:= \mu_1 - \mu_2$. Then, we take 
the difference of \eqref{uniq:11} for $i=1,2$ and test it by $u$. Noting that
$u_i(t) \in W$ for a.a.~$t\in (0,T)$ and using the fact that $\calB$ is
a maximal monotone operator from $W$ to $2^{W'}$, we have
\begin{equation}\label{uniq:13}
  \duav{\zeta_1(t)-\zeta_2(t), u_1(t) - u_2(t)}_{W',W} \ge 0 \quext{for a.a. }\,t\in (0,T),
\end{equation}
whence
\begin{equation}\label{uniq:14}
  ( \mu , u ) \ge \| A u \|^2  - ( 2\lambda - \eta ) \| \nabla u \|^2
   + \io \big( \beta(u_1)\beta'(u_1) + g(u_1) - \beta(u_2)\beta'(u_2) - g(u_2) \big) u \,\dix.
\end{equation}
Furthermore, recalling \eqref{co:33}, we infer that
\begin{equation}\no
  \ddr \big( \beta(r) \beta'(r) + g(r) \big) \ge -L \quad\perogni r\in (-1,1),
\end{equation}
where the positive constant $L$ is independent of $r$. Thus,  
\begin{align}
 & \io \big( \beta(u_1)\beta'(u_1) + g(u_1) - \beta(u_2)\beta'(u_2) - g(u_2) \big) u \,\dix \geq -L \|u\|^2\quext{for a.a. }\,t\in (0,T).
 \label{uniq:14d}
\end{align}
Since $\overline{u\zzu} = \overline{u\zzd}$, by the mass 
conservation property \eqref{st:mean}, we have $\overline{u}=0$
for a.a.  $t\in(0,T)$. Hence, recalling that the operator $A$ is invertible as it
is restricted to functions with zero-spatial mean, we
are allowed to test the difference of \eqref{E:mus} by $\mathcal{N} u$.
Noting that
\begin{equation}\label{uniq:14b}
  \duav{ A \mu, \mathcal{N} u}_{V',V}
   = \duav{ A( \mu - \overline{\mu}) , \mathcal{N} u}_{V_0',V_0}
   = \io (\mu -  \overline{\mu}) u \,\dix
   = (\mu , u),
\end{equation}
we readily obtain
\begin{equation}\label{uniq:14c}
  \frac12\ddt \| u \|_{V'}^2
   + (\mu , u) = 0.
\end{equation}
Combining \eqref{uniq:14c} with \eqref{uniq:14} and \eqref{uniq:14d}, 
we then infer
\begin{equation}\label{uniq:15}
  \frac12\ddt \| u \|_{V'}^2 + \| A u \|^2
    \le ( 2|\lambda|+|\eta |) \| \nabla u \|^2 + L \| u \|^2
    \le \frac12 \| A u \|^2 + C \| u \|_{V'}^2,
\end{equation}
where the last inequality is a consequence of the compact embeddings $W \subset V \subset V'$
and of the fact that $\| A \cdot \| + \| \cdot \|_{V'}$ is an equivalent norm on $W$.
Then, integrating \eqref{uniq:15} over $(0,T)$ and applying
Gr\"onwall's lemma, we arrive at~\eqref{cont:dep} for subdifferential solutions, which also holds for weak solutions in view of Lemma \ref{lemma:subd}. In particular, we obtain uniqueness of solutions provided that $u\zzu=u\zzd$. 

The proof of Theorem~\ref{thm:main} is now completed.

%%%%%%%%%%%%%%%%%%%%%%%%%%%%%%%%%%%%%%%%%%%%%%%%%%%%%%%%%%%%%%%%%%%%%%%%%%%%%%%%%%%%%%%%%%%%%

\section{Regularity and long-time behavior}
\label{sec:long}

In this section, we prove Theorems~\ref{thm:furth} and~\ref{thm:attra}.

%%%%%%%%%%%%%%%%%%%%%%%%%%%%%%%%%%%%%%%%%%%%%%%%%%%%%%%%%%%%%%%%%%%%%%%%%%%%%%%%%%%%%%%%%%%%%

\subsection{Proof of Theorem~\ref{thm:furth}: parabolic regularization}
\label{subssec:reg}

Our first aim is to prove parabolic regularization properties of weak solutions for strictly positive times. 

Given $h > 0$, let us introduce the difference quotient of a function $v$ by $\partial_t^h v(t)= h^{-1}[v(t+h)-v(t)]$. Applying $\partial_t^h$ to \eqref{E:mu} and testing it by $\mathcal{N} \partial_t^h u$ (noting that $\overline{\partial_t^h u}=\partial_t^h \overline{u}=0$), we get  
\begin{equation}\label{reg:11}
  \frac12\ddt \| \partial_t^h u \|_{V'}^2
   + (\partial_t^h \mu, \partial_t^h u) = 0.
\end{equation}
The second term can be computed by applying $\partial_t^h$ to \eqref{E:u}
 and testing the result by $\partial_t^h u$. This gives
\begin{align}\no
   (\partial_t^h \mu, \partial_t^h u) 
   & = \| \Delta \partial_t^h u \|^2
   + \io \partial_t^h \big( - 2\Delta \beta(u) + \beta''(u) |\nabla u|^2 \big) \partial_t^h u \,\dix\\
 \label{reg:12}
  & \quad   + \io \partial_t^h ( \beta(u)\beta'(u)+ g(u) ) \partial_t^h u \,\dix
   - (2\lambda - \eta) \| \nabla \partial_t^h u \|^2.
\end{align}
Now, exploiting as in Section~\ref{subsec:end} the convexity of the
singular diffusion part, we have 
\begin{equation}\label{reg:13}
 \io \partial_t^h \big( - 2\Delta \beta(u) + \beta''(u) |\nabla u|^2 \big) \partial_t^h u \,\dix \ge 0.
\end{equation}
Moreover, recalling \eqref{co:33}, we infer
\begin{align}
   & \io \partial_t^h (\beta(u)\beta'(u)+ g(u) ) \partial_t^h u \,\dix\nonumber\\
   &\quad = \io \int_0^1 \left[\beta'(\tau u(t+h)+(1-\tau)u(t))\right]^2 (\partial_t^h u)^2 \,\mathrm{d}\tau \,\dix\no\\
   &\qquad + \io \int_0^1 \beta(\tau u(t+h)+(1-\tau)u(t))\beta''(\tau u(t+h)+(1-\tau)u(t)) (\partial_t^h u)^2 \,\mathrm{d}\tau\, \dix\no\\
   &\qquad + \io \int_0^1  g'(\tau u(t+h)+(1-\tau)u(t)) (\partial_t^h u)^2 \,\mathrm{d}\tau\, \dix\no\\
   &\quad   \ge \io \int_0^1 \left[\beta'(\tau u(t+h)+(1-\tau)u(t))\right]^2 (\partial_t^h u)^2\, \mathrm{d}\tau \,\dix - C \|\partial_t^h u \|^2,
   \label{reg:14}
\end{align}
where $C>0$ is independent of $h$. Hence, combining \eqref{reg:11}--\eqref{reg:14}, we deduce that 
\begin{equation}\label{reg:15}
   \frac12\ddt \| \partial_t^h u \|_{V'}^2
  + \io \int_0^1 \left[\beta'(\tau u(t+h)+(1-\tau)u(t))\right]^2 (\partial_t^h u)^2 \,\mathrm{d}\tau \,\dix
  + \| \Delta \partial_t^h u  \|^2
   \le C \| \partial_t^h u  \|_V^2.
\end{equation}
Using the compact embeddings $W \subset V \subset V'$ and Ehrling's
lemma, we then have
\begin{equation}\label{reg:16}
  \ddt \| \partial_t^h u \|_{V'}^2
   + \io \int_0^1 \left[\beta'(\tau u(t+h)+(1-\tau)u(t))\right]^2 (\partial_t^h u)^2 \,\mathrm{d}\tau\, \dix
   + \| \Delta \partial_t^h u  \|^2
   \le C \| \partial_t^h u  \|^2_{V'}.
\end{equation}
Multiplying the above inequality by $t$ and integrating by parts,
we infer
\begin{align}\label{reg:17}
  & \ddt \big( t \| \partial_t^h u \|_{V'}^2 \big)
   + t \io \int_0^1 \left[\beta'(\tau u(t+h)+(1-\tau)u(t))\right]^2  (\partial_t^h u)^2 \,\mathrm{d}\tau\,\dix  
   + t \| \Delta \partial_t^h u \|^2    \no\\
  &\quad  
  \le (C t+1) \| \partial_t^h u \|^2_{V'}.
\end{align}
Hence, integrating over $(0,t)$ for $t\in(0,T)$, we obtain 
\begin{align}
& t \| \partial_t^h u \|_{V'}^2 + \int_0^t  s \io \int_0^1 \left[\beta'(\tau u(s+h)+(1-\tau)u(s))\right]^2 (\partial_t^h u)^2\, \mathrm{d}\tau \,\dix\, \mathrm{d}s+ 
\int_0^t s \| \Delta \partial_t^h u \|^2 \,\mathrm{d}s\no\\
&\quad \leq (C t+1) \int_0^t  \| \partial_t^h u(s) \|^2_{V'}\, \mathrm{d}s.
\label{reg:17b}
\end{align}
Since $\partial_t u\in L^2(0,T; V')$, it holds for a.a. $t\in (0,T)$
$$
\| \partial_t^h u(s) \|_{V'}\leq \frac{1}{h}\int_s^{s+h}\|\partial_t u(\tau)\|_{V'}\,\mathrm{d}\tau \to_{h\to 0} \|\partial_t u(s)\|_{V'}
$$
and $\|\partial_t^h u\|_{L^2(0,T;V')}\leq \|\partial_t u\|_{L^2(0,T;V')}$. 
Then, taking $h\to 0$ in \eqref{reg:17b}, in view of \eqref{reg:u} and \eqref{reg:bu1}, we deduce that for all $\tau>0$,
\begin{equation}\label{reg:18}
  \tau \| \partial_t u \|^2_{L^\infty(\tau,T;V')}
   + \tau \| \Delta \partial_t u \|_{L^2(\tau,T;H)}^2
   + \tau \| \beta(u) \|_{H^1(\tau,T;H)}^2
  \le C.
\end{equation}
This together with \eqref{E:mu} implies immediately 
\begin{align}
\tau \| \nabla \mu \|^2_{L^\infty(\tau,T;H)} \leq C,
\label{reg:18b}
\end{align}
Then the first and second conclusions in \eqref{reg:u+}
and the first conclusion in \eqref{reg:bu1+} are justified. 
Moreover,
going back to \eqref{co:22}, \eqref{co:23} and using \eqref{reg:18b},
one can easily deduce that
\begin{equation}\label{reg:19}
  \tau \| \beta(u)\beta'(u) \|_{L^1(\Omega)}
   + \tau \| \nabla\beta(u) \|^2
  \le C \tau \big( 1 + \| \nabla \mu \| \big)
  \le C,
\end{equation}
whence we deduce the second conclusion in \eqref{reg:bu1+} and also \eqref{reg:bu2+}.
Besides, recalling \eqref{st:23b}, we deduce \eqref{reg:bu3+}.
Then \eqref{st:24a} can be improved to
$\overline{\mu} \in L^\infty(\tau,T)$, which combined with \eqref{reg:18b}
yields \eqref{reg:mu+}. Finally, the third conclusion of \eqref{reg:u+} 
follows from \eqref{co:31}, \eqref{reg:18b}, \eqref{reg:19} and \eqref{reg:u}. 

Next, we prove the energy equality \eqref{diss:E2}.
To this aim, we start with observing that, thanks to the additional
regularity properties \eqref{reg:u+}--\eqref{reg:mu+}, for a.a.~$t\in (0,T)$
we can now test \eqref{E:mu} by $\mu \in V$ and test \eqref{E:u} by $u_t\in W$,
obtaining
\begin{align}
  & \| \nabla \mu \|^2
   + \duav{A^2 u, u_t}_{V',V}
   + \duav{2 A \beta(u) + \beta''(u) |\nabla u|^2, u_t}_{W',W}\nonumber\\
   &\quad 
   + \duav{\beta(u)\beta'(u)+ g(u), u_t}_{W',W}
   - (2 (\lambda - \eta) Au, u_t ) = 0.
   \label{chain:11}
\end{align}
By the standard chain rule, one has
\begin{equation}\label{chain:12}
  \duav{A^2 u, u_t}_{V',V}
   = \frac12 \ddt \| Au \|^2,
\end{equation}
whereas, recalling the chain rule formula for monotone operators in $L^2$
(see, e.g., \cite[Lemme~3.3, p.~73]{Br}) and using \eqref{reg:bu1+}, \eqref{reg:bu2+}, we readily have
\begin{align}
& \duav{\beta(u)\beta'(u)+ g(u), u_t}_{W',W}
   - (2 (\lambda - \eta) Au, u_t )\no\\
&  \quad = ( - (2 \lambda - \eta) Au, u_t ) 
          + \big( \beta(u)-\lambda u, \beta'(u) u_t \big)
          + \big( (\eta-\lambda)\beta(u)+(\lambda^2-\lambda\eta)u, u_t \big)
          \no\\
&  \quad = \ddt \io \Big( \frac12 \beta^2(u)
     - \frac12 (2\lambda - \eta) | \nabla u |^2
     + G(u) \Big) \dix,
     \label{chain:13}
\end{align}
where $G'=g$ (cf. \eqref{defig}). Finally, we deal with the more delicate term
accounting for the nonlinear diffusion part. Actually, recalling Lemma~\ref{lemma:subd},
we have, a.e.~in~$(0,T)$,
\begin{equation}\label{chain:14}
  \zeta:= 2 A \beta(u) + \beta''(u) |\nabla u|^2
   \in \calB(u),
\end{equation}
where we recall that $\calB$ represents the subdifferential of $\mathcal{J}$
with respect to the duality pairing between $W$ and $W'$.
Let us denote by $\calR :W \to W'$ the Riesz operator.
Setting $\mathcal{Z}(t):= \calR^{-1} \zeta(t)$, then we can write
\begin{equation}\label{chain:15}
  \duav{2 A \beta(u) + \beta''(u) |\nabla u|^2, u_t}_{W',W}
   = \duav{\zeta, u_t}_{W',W}
   = (\mathcal{Z}, u_t )_W
   = \ddt \mathcal{J}(u),
\end{equation}
the last equality following again from \cite[Lemme~3.3, p.~73]{Br}
applied now with respect to the scalar product $( \cdot,\cdot  )_W$
of $W$ noting that, for a.a. $t\in (0,T)$, $\mathcal{Z}(t)$ belongs to the subdifferential $\de_W \mathcal{J}(u)$ with respect to the Hilbert structure (scalar product)
of $W$. 

Collecting the above calculations \eqref{chain:11}--\eqref{chain:15}, and performing
a number of standard algebraic manipulations in order to get back
the original expression \eqref{defiE} of $\calE$, we deduce that
\begin{equation}\label{chain:16}
  \ddt \calE(u) + \| \nabla \mu \|^2
   = 0 \quext{for a.a.~}\, t\in (0,T).
\end{equation}
Since $\| \nabla \mu \|^2\in L^1(0,T)$, then integrating over $(t_1,t_2)$ for $0 < t_1 < t_2 \le T$, recalling \eqref{reg:u+}--\eqref{reg:mu+}, 
we eventually obtain the energy equality \eqref{diss:E2} for $t_1>0$, which also implies that 
$t\mapsto \calE(t)$ is absolutely continuous over
$(0,T]$. 

Here, we need to pay some more attention to the case $t_1=0$.
We notice that, if the initial datum $u_0$
is smoother, then the regularity properties \eqref{reg:u+}--\eqref{reg:mu+}
also hold for $\tau=0$, which easily implies the validity of the
energy equality starting from the initial time. 
Thus, let us consider a sequence of smooth initial data $u\zzn$
prepared in such a way that $u\zzn\to u_0$ in a suitable way 
(for instance, with respect to the distance $\mathrm{dist}_{\calX}$ as in \eqref{distXm}) and, 
in particular, $\calE(u\zzn)\to \calE(u_0)$. Let $u_n$ be the weak solution originating
from $u\zzn$. Then we have
\begin{equation}\label{chain:18}
  \calE(u_n(t)) + \int_0^t \| \nabla \mu_n(s) \|^2 \,\dis = \calE(u\zzn),\quad \forall\,t\in(0,T].
\end{equation}
Taking the $\liminf$ as $n\nearrow +\infty$ we infer
\begin{align}\label{chain:19}
  \calE(u(t)) + \int_0^t \| \nabla \mu(s) \|^2 \,\dis 
  & \le \liminf_{n\nearrow+\infty} \calE(u_n(t)) + \liminf_{n\nearrow+\infty} \int_0^t \| \nabla \mu_n(s) \|^2 \,\dis\no\\
  & \le \liminf_{n\nearrow+\infty} \left( \calE(u_n(t)) + \int_0^t \| \nabla \mu_n(s) \|^2 \,\dis\right) \no\\
  &  =\liminf_{n\nearrow+\infty} \calE(u\zzn)
   = \lim_{n\nearrow+\infty} \calE(u\zzn)
   = \calE(u_0),
\end{align}
for every $t\in(0,T]$. 

On the other hand, \eqref{reg:u} implies that $u\in C_w([0,T]; W)\cap C([0,T];H^{2-\sigma}(\Omega))$ for
every $\sigma\in (0,1/2)$. Then we can find a decreasing sequence $t_k \searrow 0$ such that 
\begin{align}
& \|\Delta u_0\|^2\leq  \liminf_{t_k\searrow 0} \|\Delta u(t_k)\|^2,\quad \|\nabla u_0\|^2=  \lim_{t_k\searrow 0} \|\nabla u(t_k)\|^2,\label{conv:1}\\
& u(t_k)\to u_0,\quad \nabla u(t_k) \to \nabla u_0,\quad \text{a.e. in}\ \Omega.\nonumber
\end{align} 
The latter implies that, a.e.~in $\Omega$, it holds 
\begin{align*}
&F(u(t_k))\to F(u_0),\quad f(u(t_k))\to f(u_0),\quad f'(u(t_k))|\nabla u(t_k)|^2\to f'(u_0)|\nabla u_0|^2.
\end{align*}
From the a.e. convergence, the boundedness of $|F(r)|$ on $[-1,1]$, the estimate $\|u(t)\|_{L^\infty(\Omega)}\leq 1$ and Lebesgue's dominated convergence theorem, we also get
\begin{align}
&\lim_{t_k\searrow 0} \io F(u(t_k))\, \dix= \io F(u_0)\, \dix.\label{conv:2}
\end{align}
Next, in view of \eqref{coerc:en}, we infer that $\|\beta(u(t_k))\|^2\leq C$, 
 and in particular from the second line of \eqref{coerc:en} with the non-negativity of $\beta'$, we also have 
$$
\io \beta'(u(t_k))|\nabla u(t_k)|^2 \dix\leq C,
$$
where these bounds may depend on $\calE(u_0)$, $\eta,\, \lambda,\, \Omega$ but are independent of $t_k$. Thus, from the a.e. convergence and Fatou's lemma, we get 
\begin{align}
&\io |f(u_0)|^2 \,\dix \leq \liminf_{t_k\searrow 0} \io |f(u(t_k))|^2 \, \dix,
\label{conv:3}\\
&\io \beta'(u_0)|\nabla u_0|^2 \, \dix \leq \liminf_{t_k\searrow 0} \io \beta'(u(t_k))|\nabla u(t_k)|^2 \,\dix,
\label{conv:4}
\end{align}
where \eqref{conv:4} together with the convergence of $\|\nabla u(t_k)\|$ in \eqref{conv:1} further yields
\begin{align}
&\io f'(u_0)|\nabla u_0|^2 \, \dix \leq \liminf_{t_k\searrow 0} \io f'(u(t_k))|\nabla u(t_k)|^2 \,\dix.\label{conv:4b}
\end{align}
Consequently, by the definition of $\calE$ and \eqref{conv:1}--\eqref{conv:4b}, we deduce that 
\begin{align}
  \calE(u_0) \le \liminf_{t_k\searrow 0} \calE(u(t_k)), \label{conv:5}
\end{align}
for the time sequence $\{t_k\}$ chosen above.
Since $\calE(u(t))$ is decreasing in time, \eqref{conv:5} actually holds for all $t\searrow 0$. Hence, from \eqref{chain:19} and \eqref{conv:5} we conclude that 
\begin{equation}\label{chain:20}
  \calE(u_0) \le \liminf_{t\searrow 0} \calE(u(t))
   \le \limsup_{t\searrow 0} \calE(u(t))
   \le \calE(u_0),
\end{equation}
whence the energy $\calE(u(t))$ is continuous at $t=0$. Therefore, the
energy equality \eqref{diss:E2} also holds for $t_1=0$, which
concludes the proof of Theorem~\ref{thm:furth}.

%%%%%%%%%%%%%%%%%%%%%%%%%%%%%%%%%%%%%%%%%%%%%%%%%%%%%%%%%%%%%%%%%%%%%%%%%%%%%%%%%%%%%%%%%%%%%

\subsection{Proof of Theorem~\ref{thm:attra}: the global attractor}
\label{subssec:attr}

Thanks to Theorem~\ref{thm:main}
and the mass-conservation property~\eqref{st:mean}, we see that weak solutions to problem \eqref{eq:mu}--\eqref{ini} generate a continuous semigroup $S(\cdot)$ on the phase-space $\calX_m$, endowed with a weaker metric given by $\mathrm{dist}_{\text{weak}}(u_1,u_2):=\|u_1-u_2\|_{V_0'}$. This implies that $S(t)$ is at least a \textit{closed semigroup} on the complete metric space $\calX_m$ in the sense of \cite{PZ07}. Therefore, the existence of a global attractor with the desired regularity follows from a standard argument in the theory of infinite-dimensional dynamical systems (see, e.g., \cite{BV,Tem}), provided that we can obtain the dissipativity and the asymptotic compactness of $S(t)$.  \smallskip

(1) \textbf{Existence of an absorbing set}. We prove the existence of a uniformly absorbing set $\mathbf{B}_0$
that is bounded in $\calX_m$. To this end, we sum the relations \eqref{st:energy}
and \eqref{co:20}. At this level we cannot take into account
the information coming from \eqref{st:11}--\eqref{st:13}. Nevertheless, using \eqref{co:21}--\eqref{co:23} (noting that the constants therein do not depend on the solution $u$), it is not difficult to arrive at
\begin{align}\no
  & \ddt \calE(u)
   + \| \Delta u \|^2
   + \| \nabla \mu \|^2
   + \kappa \| \beta(u) \beta'(u) \|_{L^1(\Omega)}
   + \kappa_1 \| \nabla \beta(u) \|^2 \\
 \no
  & \mbox{}~~~~~
   \le C \big( 1 + \| \nabla \mu \| \| \nabla u \| \big)
    + ( 2|\lambda|+ |\eta| ) \| \nabla u \|^2\\
 \no
  & \mbox{}~~~~~
   \le \frac12 \| \nabla \mu \|^2
    + C \big( 1 + \| \nabla u \|^2 \big)\\
 \no
  & \mbox{}~~~~~
   = \frac12 \| \nabla \mu \|^2
    + C \big( 1 + \| u \| \| \Delta u \| \big)\\
 \no
  & \mbox{}~~~~~
   \le \frac12 \| \nabla \mu \|^2
    + \frac12 \| \Delta u \|^2
    + C \big( 1 + \| u \|^2 \big)\\
 \label{co:61}
  & \mbox{}~~~~~
   \le \frac12 \| \nabla \mu \|^2
    + \frac12 \| \Delta u \|^2
    + C,
\end{align}
the last inequality following from the fact that $-1< u < 1$ a.e. in $\Omega$. The constant $C$ only depends on $\Omega$, $\eta$, $\lambda$, and on the initial datum $u_0$ only through the conserved quantity $m$. Thus, we end up with
\begin{equation}\label{co:62}
  \ddt \calE(u)
   + \kappa_2 \big( \| \Delta u \|^2
   + \| \nabla \mu \|^2
   + \| \beta(u) \beta'(u) \|_{L^1(\Omega)}
   + \| \nabla \beta(u) \|^2 \big)
  \le C,
\end{equation}
where $\kappa_2$ and $C$ are two uniform positive constants.

Next, recalling \eqref{defibeta} and \eqref{defia},
it is not difficult to realize that
\begin{equation}\label{co:63}
    \| \beta(u) \beta'(u) \|_{L^1(\Omega)}
    \ge \kappa_3 \| \beta(u) \|^2 - C,
\end{equation}
where $\kappa_3$, $C$ do not depend on $u$.  Then \eqref{co:62} can be more concisely rewritten
as
\begin{equation}\label{co:64}
  \ddt \calE(u)
   + \kappa_4 \calE(u)
   + \kappa_4 \| \nabla \mu \|^2
  \le C,
\end{equation}
where the positive constants $\kappa_4$ and $C$ again only depend on $\Omega$, $\eta$, $\lambda$, and on the initial datum $u_0$ only through the conserved quantity $m$.
Hence, the energy functional $\calE$ satisfies a
dissipative differential inequality, which implies that
\begin{align}
\calE(u(t))\leq \calE(u_0) e^{-\kappa_4 t}+\frac{C}{\kappa_4},\quad \forall\, t\geq 0. \label{dissi}
\end{align}
On the other hand, it is easy to verify that $\calE$ controls the ``magnitude''
of the initial data with respect
to the metric structure of $\calX_m$ both from
above and from below (cf. \eqref{coerc:en}), namely,
there exist uniform constants $\kappa_5,\,\kappa_6>0$
(we may admit their dependence on $m$, but in
fact they are independent of it, since $|m|<1$) such that
\begin{equation}\label{co:65}
  \frac14 \dist^2_{\calX} (u,0) - \kappa_6
   \le \calE(u)
   \le 4\dist^2_{\calX} (u,0) + \kappa_8,\quad \forall\, u\in\calX_m.
\end{equation}
Hence, \eqref{dissi} and \eqref{co:65} lead to the existence of a metric
bounded absorbing set $\mathbf{B}_0\subset \calX_m$. Namely, for
any bounded subset $\mathbf{B}\in \calX_m$, there exists a time $T_0=T_0(\mathbf{B})>0$ such that
$$S(t)\mathbf{B}\subset \mathbf{B}_0,\quad \forall\,t\ge T_0.$$

(2) \textbf{Asymptotic compactness of} $S(t)$. Now we prove the asymptotic compactness property of $S(t)$. The following lemma will be useful. 
\bele\label{lem:compa}
Set 
\begin{equation}\label{defiXmb}
   \widehat{\calX}_m:=\big\{ v\in W\cap H^3(\Omega):~\beta(v) \in V,~\overline{v} = m\big\},
 \end{equation}
with the distance given by 
 \begin{equation}\label{distXmb}
   \dist_{\widehat{\calX}}(v_1,v_2):= \| v_1 - v_2 \|_{H^3(\Omega)} + \| \beta(v_1) - \beta(v_2) \|_V.
 \end{equation}
 Then $\widehat{\calX}_m\subset \calX_m$ with compact immersion, namely, every closed ball in $\widehat{\calX}_m$ 
 has compact closure in $\calX_m$. 
\enle
\begin{proof}
Clearly, it is enough to show that, if $\{u_n\}$ is a bounded sequence in $\widehat{\calX}_m$ with respect to $\dist_{\widehat{\calX}}$,
then there exist a function $u\in \calX_m$
and a subsequence of $\{u_n\}$ converging to $u$ with respect to ${\rm dist}_{\calX}$. First, by the Sobolev embedding theorem, for a subsequence $u_n$ 
(not relabeled here and below), there exists a function $u\in W$, $\overline{u} = m$ such that $u_n\to u$ strongly in $W$ and a.e. in $\Omega$, which also 
implies $\beta(u_n)\to \beta(u)$ a.e.~in $\Omega$. Besides, from the uniform boundedness of $\|\beta(u_n)\|_V$, we infer that $\beta(u_n)\to \Phi$ strongly 
for some function $\Phi \in H$ and the convergence holds a.e. in $\Omega$ (possibly up to a further subsequence). Then it follows that $\beta(u)=\Phi$ a.e. in $\Omega$. 
Therefore, we have $\beta(u)\in H$ and the subsequent strong convergence 
$\beta(u_n)\to \beta(u)$ in $H$. The proof is complete.  
\end{proof}
Now, the energy identity \eqref{diss:E2}, \eqref{dissi} and \eqref{co:65} imply the uniform-in-time estimates:
\begin{align}
&\|u\|_{L^\infty(T_0,+\infty;W)}+\|\beta(u)\|_{L^\infty(T_0,+\infty;H)}\leq C,
\label{co:66}\\
&\|\partial_t u\|_{L^2(T_0,+\infty;V')}+\|\nabla \mu\|_{L^2(T_0,+\infty;H)}\leq C, 
\label{co:67} 
\end{align}
where the constant $C$ depends on $\Omega$, $\eta$, $\lambda$, and on the initial datum $u_0$ only through the conserved quantity $m$, but not on the radius of $\mathbf{B}$.
In view of the proof of Theorem~\ref{thm:furth}, below we just proceed in a formal way for simplicity. Similar to \eqref{reg:16}, we have 
\begin{equation}\label{reg:16b}
  \ddt \| \partial_t u \|_{V'}^2
   + \io \left(\beta'(u)\right)^2 (\partial_t u)^2 \dix
   + \| \Delta \partial_t u  \|^2
   \le C \| \partial_t u  \|^2_{V'},
\end{equation}
where the constant $C$ depends on $\Omega$, $\eta$, $\lambda$, and on the initial datum $u_0$ only through the conserved quantity $m$. It follows from \eqref{co:67}
that 
\begin{align}
&\int_{t}^{t+1} \|\partial_t u(s)\|_{V'}^2 \mathrm{d}s\leq C, \quad \forall\, t\geq T_0.   
\label{co:68} 
\end{align}
Therefore, the uniform Gronwall lemma (see \cite[Lemma 1.1, Chapter III]{Tem})
leads to
\begin{align}
\|\partial_t u(t)\|_{V'}^2\leq C,\quad \forall\, t\geq T_0+1. \label{co:69} 
\end{align}
By comparison, it holds 
\begin{align}
\|\nabla \mu(t)\|\leq C,\quad \forall\, t\geq T_0+1,\label{co:70} 
\end{align}
and from \eqref{co:22}, \eqref{co:23}, \eqref{co:66}, we get
\begin{align}
\|\beta(u(t))\|_{V} +\|\beta(u)\beta'(u)\|_{L^1(\Omega)} \leq C,\quad \forall\, t\geq  T_0+1. \label{co:71} 
\end{align}
The above estimate together with \eqref{co:31a}, \eqref{co:66} and \eqref{co:70} yields
\begin{align}
\|u(t)\|_{H^3(\Omega)}\leq C,\quad \forall\, t\geq  T_0+1.\label{co:72} 
\end{align}
Hence, we infer from \eqref{co:71}, \eqref{co:72} and Lemma \ref{lem:compa} the existence of a compact absorbing set $\mathbf{B}_1$ with entering time $T_1=T_0+1$. This implies the asymptotic compactness of $S(t)$ with respect to the metric for $\calX_m$ given by \eqref{distXm}. 

Thus, we can conclude that $S(t)$ admits a global attractor $\calA_m$ that is compact in $\calX_m$, with the uniform estimate \eqref{reg:attra}. The proof of Theorem~\ref{thm:attra} is complete.

%%%%%%%%%%%%%%%%%%%%%%%%%%%%%%%%%%%%%%%%%%%%%%%%%%%%%%%%%%%%%%%%%%%%%%%%%%%%%%%%%%%%%%%%%%%%%

\subsection{Proof of Proposition~\ref{prop:sepa}: the strict separation property in lower dimensions}
\label{subssec:sepa}

We note that one sufficient condition
 for the separation property \eqref{separ1} is 
\begin{align}
\|\omega\|_{L^\infty(\tau,T;L^\infty(\Omega))}\leq C,\quad \forall\,\tau\in (0,T).
\label{reg:omeLi}
\end{align}
 Indeed, in view of the arguments in \cite{MZ,GGM}, the above estimate combined with the monotonicity 
 of the singular term $\beta(u)$ enables us to show that $\|\beta(u)\|_{L^\infty(\tau,T;L^\infty(\Omega))}$ is bounded and thus \eqref{separ1} holds. Recall that we now have $\omega\in L^\infty(\tau,T;V)$  
 according to \eqref{reg:u+}--\eqref{reg:bu1+}. Then the conclusion easily follows from the Sobolev embedding theorem when the spatial dimension is one. 
 
The proof for the two dimensional case is a bit more involved. First, it follows from \eqref{reg:bu1+} that $\beta(u) \in L^\infty(\tau,T;V)$. On the other hand, thanks to the Trudinger--Moser inequaltiy in two dimensions (see e.g., \cite[Theorem~2.2]{NSY97}), we have
\begin{equation*}
	\int_{\Omega } e^{p|\beta(u)|} \mathrm{d}x \le C_{\rm TM}
	e^{C_{\rm TM} p^2 \|\beta(u)\|_{V}^2},\quad \forall\, p\in (1,\infty),
\end{equation*}
where the positive constant $C_{\rm TM}$ only depends on $\Omega$. As a consequence, using the above estimates and the simple fact $|\beta'(r)|\leq e^{2|\beta(r)|}$ for $r\in (-1,1)$, we can deduce that 
\begin{align*}
  \|\beta'(u)\|_{L^\infty(\tau,T; L^p(\Omega))} \leq C(p),\quad \forall\, p\in (1,\infty),\ \ \forall\,\tau\in (0,T),  
\end{align*}
%
%Furthermore, from the fact $|\beta''(r)|\leq 2 |\beta'(r)|^2$, $r\in (-1,1)$, we have 
%
%\begin{align*}
%  \|\beta''(u)\|_{L^\infty(\tau,T; L^p(\Omega))} \leq C(p),\quad \forall\, p\in (1,\infty),\ \ \forall\,\tau\in (0,T),  
%\end{align*}
% 
where $C(p)$ is a positive constant that may depend on the index $p$. This estimate also easily implies $f'(u)\in L^\infty(\tau,T;L^p(\Omega))$ for any $p\in(1,\infty)$. 
Now we consider the elliptic problem 
\begin{align*}
-\Delta \omega= \mu-f'(u)\omega-\eta\omega\quad \text{in}\ \ \Omega,\quad \dn \omega=0\quad \text{on}\ \ \partial\Omega.
\end{align*}
From $\mu\in L^\infty(\tau, T; V)$ (see \eqref{reg:mu+}), $\omega\in L^\infty(\tau,T;V)$ and $f'(u)\in L^\infty(\tau,T;L^p(\Omega))$, we see that $\Delta \omega\in L^\infty(\tau,T;L^p(\Omega))$ for any $p\in (1,\infty)$. Hence, by the standard elliptic regularity theory, we obtain $\omega\in L^\infty(\tau,T; W^{2,p}(\Omega))$, which together with the Sobolev embedding theorem yields the expected estimate \eqref{reg:omeLi}. 

The proof of Proposition~\ref{prop:sepa} is complete. 

\bigskip

%%%%%%%%%%%%%%%%%%%%%%%%%%%%%%%%%%%%%%%%%%%%%%%%%%%%%%%%%%%%%%%%%%%%%%%%%%%%%%%%%%%%%%%%%%%%%%%%%%%%%%%%%%%
\noindent 
\textbf{Acknowledgments.} 

G.~Schimperna has been partially supported by GNAMPA (Gruppo Nazionale per l'Analisi Matematica, la Probabilit\`a e le loro Applicazioni)
of INdAM (Istituto Nazionale di Alta Matematica). 
The research of H.~Wu is partially supported by NNSFC grant No.
11631011 and the Shanghai Center for Mathematical Sciences. 
The authors would like to thank Dr. Andrea Giorgini for helpful discussions, in particular for his suggestion on the proof of Proposition \ref{prop:sepa}.
\medskip

\noindent 
\textbf{Conflict of Interest.} 

The authors declare that they have no conflict of interest.

%%%%%%%%%%%%%%%%%%%%%%%%%%%%%%%%%%%%%%%%%%%%%%%%%%%%%%%%%%%%%%%%%%%%%%%%%%%%%%%%%%%%%%%%%%%%%%%%%%%%%%%%%%%

\end{document}